\documentclass[11pt, a4paper]{article}
\usepackage[utf8]{inputenc}
\usepackage{amssymb,amsmath,dsfont}
\usepackage{graphicx}
\usepackage{wrapfig}
\usepackage{theorem} 
\usepackage[margin=2.5cm]{geometry}
\usepackage[dvipsnames,usenames]{color}
\usepackage{epstopdf}
\usepackage[font=small,labelfont=bf]{caption}

\usepackage{float}

\usepackage[colorlinks=true, linkcolor=blue, urlcolor=magenta, citecolor=blue]{hyperref}
\newtheorem{problem}{Inverse Problem}
{\theorembodyfont{\rmfamily}\newtheorem{test}{Test}}
\newtheorem{lemma}{Lemma}[section]
\newtheorem{theorem}[lemma]{Theorem}
\newtheorem{prop}[lemma]{Proposition}
\newtheorem{coro}[lemma]{Corollary}
\newtheorem{defi}[lemma]{Definition}
{\theorembodyfont{\rmfamily}\newtheorem{remark}[lemma]{Remark}}

\newcommand{\ppp}{\partial}

\title{Reconstruction of degenerate conductivity region for parabolic equations}
\author{Piermarco Cannarsa \thanks{University of Rome Tor Vergata, Italy (e-mail: {\tt cannarsa@mat.uniroma2.it}).} \and Anna Doubova\thanks{Universidad de Sevilla, Dpto. EDAN e IMUS, Spain (e-mail: {\tt doubova@us.es}).} \and Masahiro Yamamoto\thanks{The University of Tokyo, Japan, Honorary Member of Academy of Romanian Scientists,
Correspondence member of Accademia Peloritana dei Pericolanti, Messina, Italy (e-mail:  {\tt myama@next.odn.ne.jp}).}}

\date{\today}

\begin{document}

\maketitle 

\begin{abstract}
We consider an inverse problem of reconstructing a degeneracy point in the diffusion coefficient in a one-dimensional parabolic equation by measuring the normal derivative on one side of the domain boundary. We analyze the sensitivity of the inverse problem to the initial data. We give sufficient conditions on the initial data for uniqueness and stability for the one-point measurement and show some examples of positive and negative results. On the other hand, we present more general uniqueness results, also for the identification of an initial data by measurements distributed over time. The proofs are based on an explicit form of the solution by means of Bessel functions of the first type. Finally, the theoretical results are supported by numerical experiments.
\end{abstract}

\section{Introduction}


In this paper we will consider an inverse problem of reconstruction of a degeneracy region at a point $a\in (0,1)$ for the following degenerate parabolic equation:
\begin{equation}\label{eq.cp}
\left\{\begin{array}{ll}
\partial_t w- \partial_x(|x-a| \partial_x w) = 0, & (x,t)\in(0,1)\times(0,T),  \\[1mm]
w(0,t) = 0, \quad w(1,t) = 0, &  t\in (0,T),  \\[1mm]
w(x,0) = w_0(x),  & x\in(0,1),
\end{array}\right.
\end{equation}
where $T>0$, $w_0\in C^1([0,1])$, $w_0\neq 0$,  are given. 

\medskip
Our goal is to determine or estimate the degeneracy point $a\in (0,1)$ from suitable measurements. 
Naturally, we can discuss more general cases, but we concentrate here on the one-dimensional linear equation. 

Notice that one of the difficulties in studying \eqref{eq.cp} is the fact that the energy space for such a problem in not fixed, but it depends on $a$ (see~\cite{CMV}). 

Since problem \eqref{eq.cp} is strongly degenerate, it can be decoupled into two sub-problems. More specifically, we can analyze separatly the following two problems on $(0,a)$ and $(a,1)$:
%
%
\begin{equation}\label{eq.pl}
\left\{\begin{array}{ll}
\partial_t v- \partial_x((a-x) \partial_x v) = 0, & (x,t)\in(0,a)\times(0,T),  \\[1mm]
v(0,t) = 0, \quad (a-x)\partial_x v(x,t)\big|_{x=a} = 0, &  t\in (0,T),  \\[1mm]
v(x,0) = v_0(x),  & x\in(0,a)
\end{array}\right.
\end{equation}
and 
%
%
\begin{equation}\label{eq.pr}
\left\{\begin{array}{ll}
\partial_t u- \partial_x((x-a) \partial_x u) = 0, & (x,t)\in(a,1)\times(0,T),  \\[1mm]
u(1,t) = 0, \quad (x-a)\partial_x u(x,t)\big|_{x=a} = 0, &  t\in (0,T),  \\[1mm]
u(x,0) = u_0(x),  & x\in(a,1). 
\end{array}\right.
\end{equation}

Consequently, we solve problems  \eqref{eq.pl} and \eqref{eq.pr} and get the solutions $v$ on 
$(0,a)\times(0,T)$ and $u$ on $(a,1)\times(0,T)$, respectively,  which implies that we will have the solution to \eqref{eq.cp} given by $w: = v$ for $x\in (0,a)$ and $w:=u$ for $x\in (a,1)$.

\medskip
In the sequel, we will concentrate on the analysis on the lateral problem \eqref{eq.pr}. For each $a$, let us call $u^a= u^a(x,t)$ the corresponding solution to \eqref{eq.pr}. 

A natural measurement  of the solution to \eqref{eq.pr} is  the normal derivative $\partial_xu^a(1,t)$. 
Therefore, we will consider the following inverse problem: 

\begin{problem}[Interior Degeneracy Reconstruction Problem (IRD)]
Find the degeneracy point $a\in (0,1)$ from the measurement of $\partial_x u^a(1,t)$. 
\end{problem}

\bigskip
The main issues related to this inverse problem are the following: 

\begin{description}
\item[Uniqueness:] let $u^{a_i}$, $i=1,2$ be two solutions to \eqref{eq.pr} associated to $a_i$. Assume  that the corresponding observations $\partial_x u^{a_1}(1,t)$ and $\partial_x u^{a_2}(1,t)$ coincide, i.e.
\[
\partial_x u^{a_1}(1,t) = \partial_x u^{a_2}(1,t) \quad \text{in}\quad (0,T). 
\]
Then, do we have $a_1= a_2$?
\item[Stability:]   estimate  $|a_2-a_1|$ as a function of $||\partial_x u^{a_2}(1,t)-\partial_x u^{a_1}(1,t)||_Y$ in a suitable space $Y$. 
\item[Numerical approximation:] find numerically $a$ from $\partial_x u^{a}(1,t)$. 

\end{description}

\bigskip 
Given $u_0\in C^1([0,1])$ let us introduce the map $\mu_{u_0}: [0,1)\times[0,+\infty) \mapsto \mathds{R}$ given by 
\begin{equation}\label{eq.mapmu}
\mu_{u_0}(a,t) : = \partial_x u^a(1,t),
\end{equation}
where $u^a$ satisfies \eqref{eq.pr}.

Using this map, we can reformulate our inverse problem as the follows:

\begin{problem}[Reformulation of IRD Problem]\label{IP_refor}
Does there exist an initial value $u_0\in C^1([0,1])$ such that for some constant $C>0$ and intervals $[\alpha,\beta]\subseteq [0,1)$ and $[t_0,t_1]\subseteq [0,+\infty)$ we have that 

\begin{equation}\label{eqstar}
|\mu_{u_0}(a_1,t)- \mu_{u_0}(a_2,t) | \geq C |a_1 - a_2| 
\end{equation}
for all $a_1,a_2\in [\alpha,\beta]$ and $t\in [t_0,t_1]$.

In this case, we say that $u_0$ is admissible for IRD Problem on $[\alpha,\beta]\times[t_0,t_1]$.

\end{problem}

\medskip 

Clearly, if $u_0$ is admissible for IRD Problem on $[\alpha,\beta]\times[t_0,t_1]$, then the Lipschitz stability estimate 

\begin{equation*}
 |a_1 - a_2|  \le \dfrac{1}{C}|\partial_x u^{a_2}(1,t) -\partial_x u^{a_1}(1,t)  | 
\end{equation*}
holds true for all $a_1,a_2\in [\alpha,\beta]$ and all $ t\in [t_0,t_1]$.

\medskip
Therefore, we are interested in proving the existence  of admissible initial data. Our first main goal is to analyze the sensitivity of the inverse problem to the initial data. 

\bigskip

In the last years, degenerate parabolic equations have been getting more and more attention in view of the related significant theoretical analysis and practical applications in various fields, including climatology (see Sellers~\cite{Sellers}, D\'{\i}az~\cite{diaz}, Huang~\cite{JiHuang}),  populations genetics Ethier~\cite{Ethier}, vision Citti and Manfredini~\cite{Citti}, financial mathematics~Black and Scholes~\cite{BS} and fluid dynamic Oleinik and Samokhin~\cite{Oleinik}. See also Cannarsa, Martinez and Vancostenoble~\cite{CMV} and the references therein.

\medskip
Inverse problems are the kind of problems referred to as ill-posed in the Hadamard sense (see \cite{Hadamard}). That is, the solution either does not exist, or is non-unique and/or small errors in the given data can lead to large errors in the calculated solutions. 

\medskip
The related literature concerning inverse problems for degenerate parabolic PDEs, in contrast to its essential relevance and practical applicability, is rather scarce and recent. 

For example, the inverse  source problem was considered in Tort~\cite{Tort1}, Cannarsa, Tort, and Yamamoto~\cite{CTY}, Cannarsa, Martinez and Vancostenoble~\cite{CMV}, Deng et al.~\cite{Deng} and Hussein et al.~\cite{HusseinLesnic}, where numerical reconstruction was also considered.  The inverse problem of recovering the first-order coefficient of degenerate parabolic equations is analyzed in~Deng and Yang \cite{Deng2} and Kamynin~\cite{Kamynin}. An inverse diffusion problem was considered in~\cite{Tort2}, where a constant diffusion coefficient was recovered from measurement of second order derivatives by means of Carleman estimates in Cannarsa, Tort and Yamamoto~\cite{CTY}. In a recent paper~\cite{CDY} Cannarsa, Doubova, Yamamoto have analyzed several inverse problems of reconstruction of degenerate diffusion coefficient in a parabolic equation.  

\medskip
Regarding the inverse problems of determination of the spatially varying coefficients like conductivity and source terms,  there has been a substantial amount of work for non-degenerate parabolic equations. Just to illustrate some of them, we refer  to the chapters by Isakov~\cite{Is}, and Yamamoto~\cite{Ya}. See also the references therein.  Particularly, for the same type of inverse problems for one-dimensional regular parabolic equations, we refer to Murayama~\cite{Mura},  and
Pierce \cite{Pier}, Suzuki and Murayama \cite{SuMu}. 

\medskip 
In addition, most of the techniques which have been devised to deal with non-degenerate parabolic equations are in general not applicable in the degenerate case. For example, for strongly degenerate operators, the trace of the cornormal derivative must vanish in the part of the boundary where the ellipticity fails, and so no useful measurements are obtained. 

\medskip

On the other hand, the second main issue of this paper is the reconstruction of degeneracy. The objective is to compute approximations of the degeneracy point as well as solutions to the inverse problems using observation data. 

An effective approach to accomplish this, as illustrated here below, is to reformulate the search for the degeneracy as an extremal problem. This is nowadays classical and has been widely implemented in many situations; see for instance Lavrentiev et al.~\cite{Lavrentiev}, Samarskii and  Vabishchevich~\cite{Samarskii}, Vogel~\cite{Vogel} and Cannarsa, Doubova, Yamamoto~\cite{CDY}. 

 \bigskip
 The paper is organized  as follows. In Section~\ref{sec:weelposedness}, we consider the well-posedness of the corresponding forward problem. In Section~\ref{sec:eigenvalues}, we present the main result giving an expression of the normal derivative obtained by performing explicit computations that use the Bessel function of the first kind.  In Section~\ref{sec.stablility} we establish a  Lipschitz stability result with one point measurement. 
 Section~\ref{sec.examples} will be devoted to the presentation of some examples of initial data for which we can have a stability estimate making explicit the stability constant. In addition, we will see that there exist initial data for which we can obtain Lipschitz stability estimate for an arbitrary small time. In Section~\ref{sec:uniqueness} we will present general uniqueness results for distributed  measurements over the whole  time interval. Finally, in Section~\ref{sec:numerics}, in order to illustrate the theoretical results obtained in the previous sections, we perform some satisfactory numerical experiments corresponding to the considered inverse problem.  


 \section{Well-posedness }\label{sec:weelposedness}

We start recalling the natural functions spaces where the  problem can be set. Let us call $H=L^2(0,1)$. For all $a\in (0,1)$ we consider the following function spaces: 
\begin{equation}\label{eq.spaces}
H_a^1(0,1) = \Big\{ u\in H \,  |\,    \displaystyle \int_0^1 (x-a)|\partial_x u|^2\, dx<\infty,\,   u(1) =  0 \Big\}
\end{equation}
and 
\begin{equation*}
H_a^2(0,1) =
 \{ u\in H^1_a(0,1) \,  |\,  x\mapsto (x-a) \partial_x u\in H_a^1(0,1)  \}.
\end{equation*}

\begin{remark}[Neumann boundary condition]\label{Neumann boundary condition}\label{rem:Neumann}
For $u\in H^2_a(0,1)$ we have $(x-a) \partial_x u|_{x=a}=0$. Indeed, if
$(x-a) \partial_x u(x) \to L$ when $x\to a$, then $(x-a)|\partial_x u(x) |^2 \sim L^2/(x-a)$ and, therefore $L=0$ otherwise  $u\notin H^1_a(0,1)$.
\end{remark}

Problem \eqref{eq.cp}  can be recast  in the abstract form 
\begin{equation}
\label{eq:HCP}
\begin{cases}
 u'(t)=Au(t) & t\ge 0,
 \vspace{.1cm}
 \\
 u(0)=u_0
\end{cases}
\end{equation}
by introducing the linear operator $A:D(A)\subset H \to H$ defined by
\begin{equation}
\label{eq:A}
D(A) = H_a^2(0,1), \qquad Au = \partial_x((x-a) \partial_x u)\quad \text{for } u\in D(A).
\end{equation}
Therefore, we can rely on~\cite{MarVan} to obtain a well-posedness result for  \eqref{eq:HCP}. 
Hence,  one can prove the following, where we recall that $D((-A)^{1/2})=H_a^1(0,1)$ for the operator $A$ defined by \eqref{eq:A}.

\begin{prop}\label{th2.3}
The operator $A$ is an infinitesimal generator of a strongly continuous semigroup of contractions, $e^{tA}$. Moreover, $e^{tA}$ is analytic.
Therefore, for any $u_0\in H$,  the solution $u(x,t) = (e^{tA}u_0)(x)$ of problem  \eqref{eq:HCP} satisfies:

\begin{itemize}
\item[(i)] $u\in C^0([0,\infty);H)$,
\item[(ii)] $t\mapsto u(\cdot, t)$ is analytic as a map $(0,+\infty) \longrightarrow H$,
\item[(iii)] $u(t) \in \bigcap_{n\geq 1} D(A^n)$ for all  $ t>0$. 
\end{itemize}
Furthermore, if $u_0\in D((-A)^{1/2})$, then the mild solution
\begin{equation*}
u(t):=e^{tA}u_0 \qquad t\in[0,T] 
\end{equation*}
of problem \eqref{eq:HCP} belongs to $H^1(0,T;H)\cap C([0,T];D((-A)^{1/2}))\cap L^2(0,T;D(A))$ and satisfies the equation \eqref{eq:HCP} for a.e. $t\in[0,T]$.
\end{prop}

Consequently, the equation in \eqref{eq.pr} is satisfied in classical sense on $(a,1)\times(0,+\infty)$, as well as boundary conditions, taking into account  Remark~\ref{Neumann boundary condition}. As for the initial condition, we recall that $u_0\in L^2(0,1)$ implies  $u\in C^0([0,+\infty);L^2(0,1))$ and 
$u_0\in H^1_a(0,1)$ implies  $u\in C^0([0,+\infty);H^1_a(0,1))$. Also, notice that the function $t\mapsto \partial_x u^a(1,t)$ is analytic for all $t>0$, since $u$ is analytic for all $t>0$. Hereafter, we will assume solutions as smooth as required.

 \section{Computation of the normal derivative}\label{sec:eigenvalues}
 
 In this section we will perform the explicit computation of the normal derivative $\partial_xu^a(1,t)$, where $u^a(x,t)$ is the solution to \eqref{eq.pr}. In order to do this we will need an explicit expression of solutions to the associated eigenfunctions and eigenvalues in terms of the Bessel functions of the first kind.

Let us indicate here also the dependence on the initial data of the solution to \eqref{eq.pr} by putting  $u^a: = u^a(x,t;u_0)$.  For all $x\in (a,1)$, let us perform the change of variables 

 \begin{equation}\label{eq.change}
y = \dfrac{x-a}{1-a},\quad  x = a + (1-a)y
 \end{equation}
and
  \begin{equation*}
u^a(x,t;u_0)  = v\Big(\dfrac{x-a}{1-a}, t; v_0^a\Big), \quad x\in (a,1),
 \end{equation*}
where
 \begin{equation*}
v_0^a(y)  = u_0(a + (1-a)y), \quad y\in (0,1).
 \end{equation*}
We have 

\[
\partial_x u^a(1,t;u_0) = \dfrac{1}{1-a}\partial_y v(1,t;v_0^a).
\]

\noindent
Therefore, $v = v(y,t;v_0^a)$ satisfies 
\begin{equation}\label{eq.prv}
\left\{\begin{array}{ll}
\partial_t v- \dfrac{1}{1-a}\partial_y(y\partial_y v) = 0, & (y,t)\in(0,1)\times(0,T),  \\[1mm]
v(1,t) = 0, \quad y\partial_y v(y,t)\big|_{x=0} = 0, &  t\in (0,T),  \\[1mm]
v(y,0) = v_0^a(y),  & y\in(0,1). 
\end{array}\right.
\end{equation}

Henceforth we set 

\[
\phi'(x) = \dfrac{d \phi}{dx}(x), \quad \text{etc.} 
\]
 
 \medskip
 Let us consider the eigenvalue problem: 

 \begin{equation}\label{pb.autofun}
 \begin{cases}
 - ((x-a) \phi'(x))' = \lambda \phi,  \quad x\in (a,1),
 \\[2mm]
 \phi(1) = 0,\quad   (x-a)\phi'(x)\big|_{x=a} = 0.
 \end{cases}
 \end{equation}
We transform the problem \eqref{pb.autofun} into a problem on $(0,1)$ by the change of variables \eqref{eq.change} and we obtain that $\psi(y) = \phi(a+(1-a)y)$ satisfies
 
  \begin{equation}\label{pb.autofunch}
 \begin{cases}
 - (y \psi'(x))' = \lambda (1-a) \psi,  \quad y\in (0,1),
 \\[2mm]
 \psi(1) = 0,\quad   y\psi'(x)\big|_{y=0} = 0.
 \end{cases}
 \end{equation}
 
\begin{figure}[h!]
 \centering 
 \includegraphics[width = 8cm]{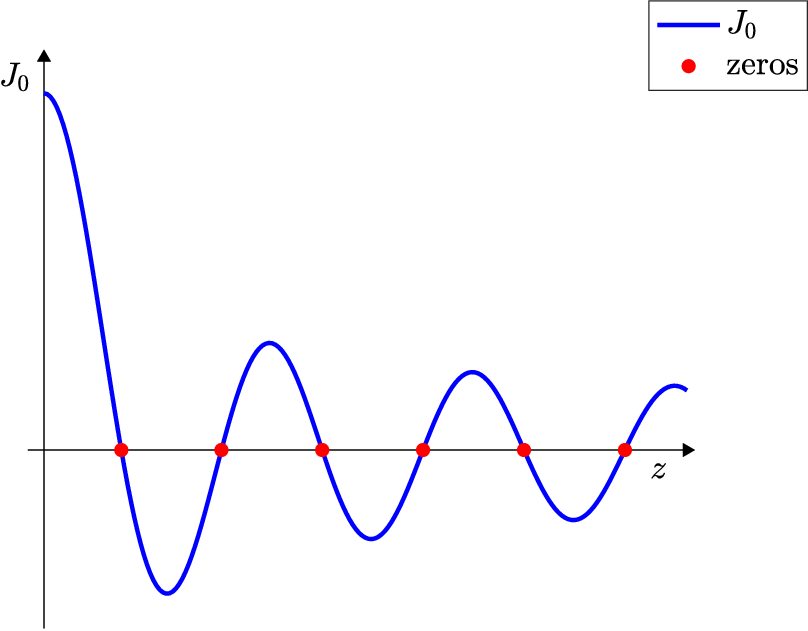}
 \caption{Bessel function $J_0$ of the first kind.}
 \label{fig.BesselJ0}
 \end{figure}

The structure of the eigenvalues $\lambda_n$ and the eigenfunctions $\psi_n$ for the problem \eqref{pb.autofun} is described by the Bessel function  $J_0$ of the first kind (see~\cite{Lebedev}): 
 \begin{equation}\label{eq.besselJ0}
J_0(z) = \displaystyle\sum_{k=0}^\infty \dfrac{(-1)^k}{(k!)^2}\Big( \frac{z}{2}\Big)^{2k}, \quad z\ge 0,
 \end{equation}
 which is a solution of the Bessel equation 
 \[
 z^2 J_0''(z) + z J_0'(z) + z^2 J_0(z) =0, \quad z\ge 0. 
 \]
 (see Figure~\ref{fig.BesselJ0}). 
 
 Notice that  the Bessel functions $J_n$, for all $n = 0,1,2, \dots$ are expressed 
 \begin{equation}\label{eq.besselJn}
J_n(z) = \displaystyle\sum_{k=0}^\infty \dfrac{(-1)^k}{k! (n+k)!}\Big( \frac{z}{2}\Big)^{n+2k} \quad n\ge 0
 \end{equation}
 
\begin{figure}[h!]
 \centering 
 \includegraphics[width = 8cm]{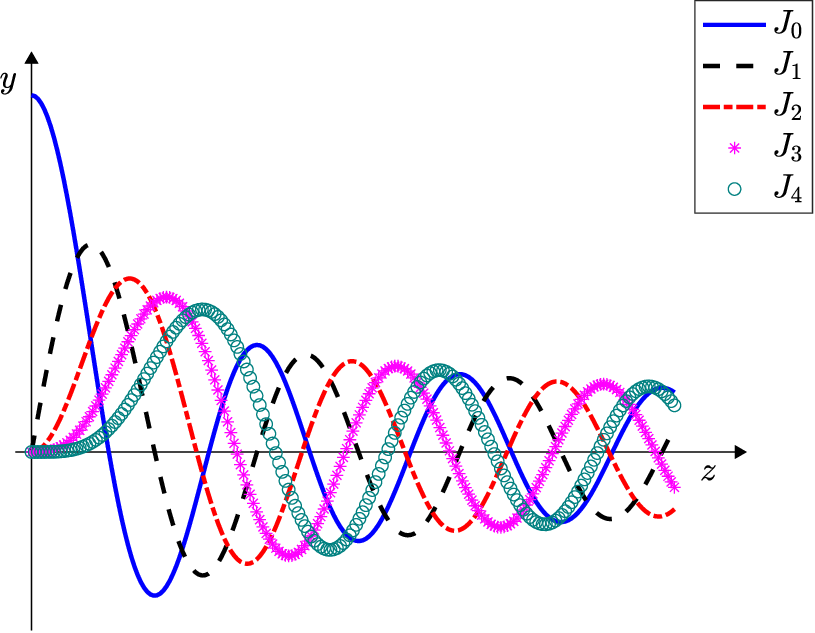}
 \caption{Bessel functions $J_n$ of the first kind.}
 \label{fig.Bessel}
 \end{figure}

 \noindent 
 (see Figure~\ref{fig.Bessel} and Lemma~\ref{lemma.Bessel}). 
 
 \medskip
 
 The following properties of the Bessel functions will be used later:
 
\begin{lemma}[Properties of Bessel functions]\label{lemma.Bessel}
Let $J_n(z)$, $n\ge 0$ be the Bessel functions of the first kind given by  \eqref{eq.besselJn} and  let us denote by $\{j_n\}_{n\ge 1}$ the sequence of positive increasing zeros of $J_0$ (see Figure~\ref{fig.BesselJ0}), i.e. $J_0(j_n)=0$, with $0<j_1<j_2<\cdots$. 

Then, the following holds: 
\begin{enumerate}
\item[\textit{a)}] $ \dfrac{d}{dz}\big (z^n J_n(z)\big) =  z^n J_{n-1}(z)$, $n = 1, 2,\dots$.
\item[\textit{b)}] $J_{n-1}(z) + J_{n+1}(z) =\dfrac{2n}{z} J_n(z)$ and $J_{n-1}(z) - J_{n+1}(z) =2J'_n(z)$, $n = 1, 2,\dots$ 
\item[\textit{c)}] $J_0'(z) = -J_1(z)$.
\item[\textit{d)}]  $\displaystyle\int_0^{j_n} s J_0^2(s) \, ds = \dfrac{j_n^2}{2} J_1^2(j_n) = \dfrac{j_n^2}{2} J_0'(j_n)^2$.
\item[\textit{e)}] $\displaystyle\int_0^{x} s J_0(s) \, ds = x J_1(x) = -x J_0'(x)$.
\item[\textit{f)}] The following bounds on the zeros of $J_0$ hold: 
$\pi\Big(n -\dfrac{1}{4}\Big) \le j_n \le \pi\Big(n -\dfrac{1}{8}\Big)$ for all $n\ge 1$.
\item[\textit{g)}]  $\displaystyle\sup_{n\in \mathds{N}} \vert J_0'(j_n) \vert <+ \infty$.
\end{enumerate} 

\end{lemma}

\medskip
\noindent
The proof of Lemma~\ref{lemma.Bessel} is given in Appendix~\ref{sec.appendix}.

\bigskip 
 
 The main result of this section is the following. 
\begin{theorem}\label{prop.normder}
Assume $u_0\in L^2(0,1)$ and that $\{j_n\}_{n\ge 1}$ is the sequence of positive zeros of $J_0$. Then, the following holds: 

\begin{enumerate}
\item[a)]  The solution to $u^a$ to \eqref{eq.pr} is given by the following: 
\begin{equation}\label{eq.solua}
 u^a(x,t) = \displaystyle\sum_{n = 1}^\infty\dfrac{2(1-a)}{j_n^2 (J_0'(j_n))^2} f_n(t,a) \, U_n^0(a) \, J_0 \Big(j_n \sqrt{\dfrac{x-a}{1-a}}\Big),
\end{equation}
where 
 \begin{equation}\label{eq.fna}
f_n(t,a) : = e^{-\big( \frac{j_n}{2}\big)^2\frac{t}{1-a}}
 \end{equation}
and 
 \begin{equation}\label{eq.U0n}
U^0_n(a) : =  \displaystyle \int_{0}^{j_n} u_0\Big(a+(1-a)\frac{s^2}{j_n^2} \Big)s J_0 (s) \, ds. 
 \end{equation}

\item[b)] As a consequence of the previous point, we have

 \begin{equation}\label{eq.normderuxa}
\partial_x u^a(1,t)  =  \displaystyle\sum_{n = 1}^\infty \dfrac{f_n(t,a)}{j_n J'_0(j_n)}U_n^0(a).
 \end{equation} 
\end{enumerate} 
 \end{theorem}
  
 Notice that where the map $a\mapsto f_n(t,a)$ is strictly decreasing for all $t>0$ and $n\ge 1$. 
 
 \bigskip
 \noindent 
 \textbf{Proof of Theorem~\ref{prop.normder}:}
\textit{a):} It is known (see~\cite{Gueye} and~\cite{CMV2}) that the solution to \eqref{pb.autofunch} is given by the eigenfunctions and the eigenvalues of the following form:
 \begin{equation}\label{eq.autofunB}
(1-a) \lambda_n = \Big(\dfrac{j_n}{2}\Big)^2, \quad  \psi_n(y) = c_n J_0(j_n \sqrt{y}),
 \end{equation}
with

\[
\dfrac{1}{c_n^2} = \int_0^1 J_0^2(j_n \sqrt{y}) \, dy = \dfrac{2}{j_n^2} \int_0^{j_n} s J_0^2(s)\, ds,
\]

\noindent
where we have made the change of variable $s = j_n \sqrt{y}$ in the integral. Taking into account Lemma~\ref{lemma.Bessel}, we have
\[
\frac{1}{c_n^2} = J_0'(j_n)^2, \quad \text{that is,} \quad c_n^2 = \dfrac{1}{J'_0(j_n)^2} = \dfrac{1}{J'_1(j_n)^2}.
\]
Therefore, returning to \eqref{pb.autofun} we obtain that 

 \begin{equation}\label{eq.autofunB2}
\lambda_n = \dfrac{1}{1-a}\Big(\dfrac{j_n}{2}\Big)^2, 
\quad  
\phi_n(x) = \dfrac{1}{|J_0'(j_n)|}J_0 \Big(j_n \sqrt{\dfrac{x-a}{1-a}}\Big).
 \end{equation}
 
\bigskip

We can represent now the solution $u^a$ to \eqref{eq.pr} as follows: 
 \begin{equation}\label{eq.solpr}
 u^a(x,t) = \displaystyle \sum_{n=1}^\infty e^{-\lambda_n t } u^0_n \phi_n(x), 
 \end{equation}
where $u^0_n$ is given by
\[
u^0_n = \displaystyle \int_{a}^1 u_0(x) \phi_n(x) \, dx.  
\]
Taking into account \eqref{eq.autofunB2} and performing the change of variables 
$s = j_n\sqrt{\frac{x-a}{1-a}}$, we obtain
 
 \begin{equation}\label{eq.uon}
\begin{split}
u^0_n = \displaystyle \int_{a}^1 u_0(x) \phi_n(x) \, dx 
 & = 
 \dfrac{1}{|J_0'(j_n)|}  \displaystyle \int_{a}^1 u_0(x) J_0 \Big(j_n \sqrt{\dfrac{x-a}{1-a}}\Big) \, dx  
 \\[3mm]
 & =  \dfrac{1}{|J_0'(j_n)|}  \displaystyle \int_{0}^{j_n} u_0\Big(a+(1-a)\frac{s^2}{j_n^2} \Big)J_0 (s) 2(1-a) \dfrac{s}{j_n^2} \, ds
  \\[3mm]
 & : =
 \dfrac{2(1-a)}{j_n^2|J_0'(j_n)|}  U_n^0(a), 
 \end{split} 
\end{equation}
 
 where $U_n^a$ is given by \eqref{eq.U0n}. From \eqref{eq.solpr}, using  \eqref{eq.autofunB2} and \eqref{eq.uon}, we deduce \eqref{eq.solua}. 

 \medskip
 \noindent 
 \textit{b)} Taking into account \eqref{eq.solua}, we obtain
 \[
  \partial_x u^a(x,t) = \displaystyle\sum_{n = 1}^\infty
  \dfrac{f_n(t,a)}{j_n (J_0'(j_n))^2} \, U_n^0(a) \, J'_0 \Big(j_n \sqrt{\dfrac{x-a}{1-a}}\Big)\Big(\dfrac{x-a}{1-a}\Big)^{-1/2}.
 \]
Hence, evaluating for $x = 1$, we easily obtain \eqref{eq.normderuxa}. 
 
Finally, since 
\[
\partial_a f_n(t,a) =  - \Big( \dfrac{j_n}{2}\Big)^2 \dfrac{t}{(1-a)^2} f_n(t,a) <0,
\]
 we also deduce that the map $a\mapsto f_n(t,a)$ is strictly decreasing for all $t>0$ and $n\ge 1$.
 
 \hfill$\blacksquare$
 
 \bigskip
 \section{Lipschitz stability with one point measurement  }\label{sec.stablility}
 
 Based on the explicit expression of the solution given in Theorem~\ref{prop.normder} we will present a stability result for one point measurement. We have the following. 
 
\begin{theorem}\label{th.stabillity}
Assume $u_0\in C^0([0,1])$. Let   $u^{a_1}$ and $u^{a_2}$ be the solutions to \eqref{eq.pr} corresponding to degeneracy points $a_1$ and $a_2$, respectively. Assume that there exist $\delta>0$ and $\beta\in [0,1)$ such that 

\begin{equation}\label{th2.assump}
|U_1^0(a)|\ge \delta \quad \forall\, a\in [0,\beta],
\end{equation}
with $U_1^0(a)$ given by \eqref{eq.U0n}. Then, there exist $T(u_0,\beta)>0$ and a constant $C>0$ such that for all  $t\in  [T(u_0), t_1]$, $t_1>T(u_0,\beta)$ the following stability estimate holds:

  \begin{equation}\label{eq.th} 
 |a_2-a_1| \le   C |\partial_x u^{a_2}(1,t) - \partial_x u^{a_1}(1,t)| \quad \forall\, a_1,a_2\in [0,\beta].
 \end{equation}
\end{theorem}

\medskip
 \begin{remark}
Notice that assumption \eqref{th2.assump} is satisfied for any $\beta\in [0,1)$ if $|u_0|>0$ for all $x\in (0,1)$. Indeed, assume $|u_0|>0$ for all $x\in (0,1)$. then $U^0_1(a) : =  \displaystyle \int_{0}^{j_1} u_0\Big(a+(1-a)\frac{s^2}{j_1^2} \Big)s J_0 (s) \, ds$ cannot approach zeros if $a\in [0,\beta]$ because the integrand is either strictly positive or strictly negative on $(0,j_1)$.

 \end{remark}
 
 \medskip
 \noindent
 \textbf{Proof of Theorem~\ref{th.stabillity}:}  Let assume to fix the idea that $U_1^0(a)\ge \delta \quad \forall\, a\in [0,\beta]$. Using the explicit formula \eqref{eq.normderuxa}, we compute 
 
 \begin{equation}\label{eq1.th2}
 \begin{split}
\partial_a(\partial_x u^a(1,t)) & =  \displaystyle\sum_{n = 1}^\infty \dfrac{e^{-\lambda_n t}}{j_n J'_0(j_n)}
\Big[(U_n^0)'(a) - \dfrac{j_n^2 t}{4(1-a)^2} U_n^0(a)\Big] 
\\[3mm]
& =  \dfrac{e^{-\lambda_1 t}}{j_1 J'_0(j_1)}
\Big[(U_1^0)'(a) - \dfrac{j_1^2 t}{4(1-a)^2} U_1^0(a)\Big] 
\\[3mm]
&  \quad+ \displaystyle\sum_{n = 2}^\infty \dfrac{e^{-\lambda_n t}}{j_n J'_0(j_n)}
\Big[(U_n^0)'(a) - \dfrac{j_n^2 t}{4(1-a)^2} U_n^0(a)\Big],
\end{split}
 \end{equation}
 where  $U_n^0(a)$ is given by \eqref{eq.U0n} and $\lambda_n$ by \eqref{eq.autofunB2}, i.e.,
 \[
 \lambda_n = \dfrac{j_n^2}{4(1-a)}, \quad 
U^0_n(a) =  \displaystyle \int_{0}^{j_n} u_0\Big(a+(1-a)\frac{s^2}{j_n^2} \Big)s J_0 (s) \, ds. 
 \]
 
 \noindent
 Let us estimate the first term in \eqref{eq1.th2}. We have
 
  \begin{equation}\label{eq2.th2}
 \begin{split}
(U_1^0)'(a) & - \dfrac{j_1^2 t}{4(1-a)^2} U_1^0(a)
\\[3mm]
&  =
\displaystyle \int_{0}^{j_1} u'_0\Big(a+(1-a)\frac{s^2}{j_1^2} \Big)\Big(1-\dfrac{s^2}{j_1^2} \Big)s J_0 (s) \, ds
-\dfrac{j_1^2 t}{4(1-a)^2} U_1^0(a) 
\\[3mm]
& \le ||u'_0||_\infty\displaystyle \int_{0}^{j_1} s J_0(s)\, ds - \dfrac{j_1^2 \delta}{4(1-\beta)^2} t
\\[3mm]
& = -j_1 J'_0(j_1)  ||u'_0||_\infty -\dfrac{j_1^2 \delta}{4(1-\beta)^2} t \quad \forall\, a\in [0,\beta].
\end{split}
 \end{equation}
 Now, observe that for all $L>0$ we have 
 
  \begin{equation}\label{eq3.th2}
 (U_1^0)'(a) - \dfrac{j_1^2 t}{4(1-a)^2} U_1^0(a) \le -L\quad \forall\, a\in [0,\beta], \quad
 \forall\, t\ge \bar{t}, 
 \end{equation}
 where $\bar{t}=t_L(u_0,\beta,L)>0$ is given by 
   \begin{equation}\label{eqTu0.th2}
\bar{t} =  \dfrac{4(1-\beta)^2}{j_1^2 \delta} \Big( L - j_1 J'_0(j_1)||u'_0||_\infty  \Big).
 \end{equation}
\noindent 
 Note that $J'_0(j_1)<0$.
 
Therefore, taking into account \eqref{eq2.th2} and \eqref{eq3.th2} in \eqref{eq1.th2}, we obtain  
\begin{equation}\label{eq4.th2}
 \begin{split}
\partial_a(\partial_x u^a(1,t)) & \ge e^{-\lambda_1 t} \Bigg[ \dfrac{L}{|j_1 J_0'(j_1)|} 
-
 \displaystyle\sum_{n = 2}^\infty \dfrac{e^{-(\lambda_n -\lambda_1)t}}{|j_n J'_0(j_n)|}
\Big|(U_n^0)'(a) - \dfrac{j_n^2 t}{4(1-a)^2} U_n^0(a)  \Big|\Bigg]
\end{split}
 \end{equation}
 for all  $t\ge \bar{t}$.
 
 \medskip
 \noindent
 Therefore, we have to show that the second term on the right-hand side of \eqref{eq4.th2} is small for $t$ large enough. We have that
 
 \begin{equation}\label{eq5.th2}
 \begin{split}
\Big |(U_n^0)'(a) & - \dfrac{j_n^2 t}{4(1-a)^2} U_n^0(a)  \Big| 
\\[3mm]
& = \Bigg |\int_{0}^{j_n} \Bigg [ u'_0\Big(a+(1-a)\frac{s^2}{j_n^2} \Big)\Big(1-\dfrac{s^2}{j_n^2} \Big)-
 \dfrac{j_n^2 t}{4(1-a)^2}u_0\Big(a+(1-a)\frac{s^2}{j_n^2}\Big) \Bigg] s J_0 (s) \, ds \Bigg|
 \\[3mm]
& \leq \Big( ||u'_0||_\infty +  \dfrac{j_n^2 t}{4(1-a)^2} ||u_0||_\infty  \Big)  \int_{0}^{j_n}s |J_0 (s)| \, ds
 \\[3mm]
& \leq \dfrac{j_n^2}{2}K_n(t),
\end{split}
 \end{equation}
where we have taken into account that $|J_0 (s)|\le 1$ and set
\[
K_n(t) = ||u'_0||_\infty +  \dfrac{j_n^2 t}{4(1-a)^2} ||u_0||_\infty.
\]
Hence, the second term in the brackets in \eqref{eq4.th2} can be estimated as follows:

 \begin{equation}\label{eq6.th2}
  \displaystyle\sum_{n = 2}^\infty \dfrac{e^{-(\lambda_n -\lambda_1)t}}{j_n |J'_0(j_n)|}
\Big|(U_n^0)'(a) - \dfrac{j_n^2 t}{4(1-a)^2} U_n^0(a) \Big | 
 \\[3mm]
 \le
 \displaystyle\sum_{n = 2}^\infty  \dfrac{1}{j_n |J'_0(j_n)|} \dfrac{j_n^2}{2}\, K_n(t) \, e^{-(\lambda_n -\lambda_1)t} : = R.
 \end{equation}
  
 Since by Lemma~\ref{lemma.Bessel} d) we have  $\int_{0}^{j_n}s J^2_0 (s) \, ds = \frac{1}{2}j_n J'_0(j_n)^2$  we conclude that 
 \[
 \lim_{n\to \infty } j_n |J'_0(j_n)| \ge M >0.
 \]
 Therefore, we can estimate the expression of $R$ in \eqref{eq6.th2} in the following way:
 
 \begin{equation}\label{eq7.th2}
 \begin{split}
 R & \le \dfrac{1}{2M} \displaystyle\sum_{n = 2}^\infty  j_n^2 \, K_n(t) \, e^{-(\lambda_n -\lambda_1)t} 
 \\[3mm]
&  \le
\dfrac{1}{2M} \Big[   ||u'_0||_\infty \displaystyle\sum_{n = 2}^\infty  j_n^2 e^{-(\lambda_n -\lambda_1)t}   
+
 \dfrac{t  ||u_0||_\infty}{4(1-a)^2} \displaystyle\sum_{n = 2}^\infty  j_n^4 e^{-(\lambda_n -\lambda_1)t}   \Big].
 \end{split}
 \end{equation}
 \noindent 
Moreover, 

\[
j_n^2 e^{-\lambda_n t} \le e^{j_n^2} e^{-\frac{t j_n^2}{4(1-a)}} \le e^{\frac{4(1-a)-t}{4(1-a)}j_n^2}
\le 
e^{-\frac{t j_n^2}{8(1-a)}}\quad \forall\, t\ge 8(1-a)
\]  
and (since $x^2\le e^x$ $\forall\, x\ge 0$) we have 
\[
j_n^4 e^{-\lambda_n t} \le e^{j_n^2} e^{-\frac{t j_n^2}{4(1-a)}} \le 
e^{-\frac{t j_n^2}{8(1-a)}}\quad \forall\, t\ge 8(1-a).
\]
\noindent
In conclusion, considering these estimates in \eqref{eq7.th2} we reach

 \begin{equation}\label{eq8.th2}
 \begin{split}
 R &  \le
\dfrac{1}{2M}  e^{\lambda_1 t}\Big[   ||u'_0||_\infty    
+
 \dfrac{t  ||u_0||_\infty}{4(1-a)^2}    \Big]
 \displaystyle\sum_{n = 2}^\infty   e^{\frac{-tj_n^2}{8(1-a)}}
 \\[3mm]
 &  \le
 \dfrac{1}{2M}  \Big[   ||u'_0||_\infty    
+
 \dfrac{t  ||u_0||_\infty}{4(1-a)^2}    \Big]
 \displaystyle\sum_{n = 2}^\infty  e^{\big(\lambda_1-\frac{j_n^2}{8(1-a)}\big)t}  \quad  \forall\,  t\ge 8(1-a).
 \end{split}
 \end{equation}
\noindent 
We claim that

 \begin{equation}\label{eq9.th2}
\lim_{t\to +\infty} \dfrac{1}{2M}  \Big[   ||u'_0||_\infty    
+
 \dfrac{t  ||u_0||_\infty}{4(1-a)^2}    \Big]
 \displaystyle\sum_{n = 2}^\infty  e^{\big(\lambda_1-\frac{j_n^2}{8(1-a)}\big)t} = 0. 
 \end{equation}
 \noindent 
Indeed, by the monotone convergence theorem we have that
  
   \begin{equation}\label{eq10.th2}
 \displaystyle\sum_{n = 2}^\infty  e^{\big(\lambda_1-\frac{j_n^2}{8(1-a)}\big)t} 
 \le 
  \displaystyle\sum_{n = 2}^\infty  e^{\big(\lambda_1-\frac{j_n^2}{8}\big)t} \to 0,  \quad t\to \infty.
 \end{equation}
As for the second term in \eqref{eq9.th2}, we have that (for $n$ sufficiently large so that 
$\lambda_1-\frac{j_n^2}{8(1-a)}\le -\frac{j_n^2}{16(1-a)}$)

  \begin{equation}\label{eq11.th2}
t e^{\big(\lambda_1-\frac{j_n^2}{8(1-a)}\big)t} 
 \le 
 t e^{-\frac{j_n^2 t}{16(1-a)}} \dfrac{j_n^2}{16(1-a)} \frac{16(1-a)}{j_n^2} \le \frac{16(1-a)}{e} \dfrac{1}{j_n^2}
 \end{equation} 
 which is summable for $j_n^2\sim n^2$. Notice that in \eqref{eq11.th2} we have used that $ t e^{-\frac{j_n^2 t}{16(1-a)}} \frac{j_n^2}{16(1-a)} \le e^{-1}$, since 
 $xe^{-x} \le e^{-1}$.
 
 Therefore,  Lebesgue's dominant convergence theorem yields
 \[
\lim_{t\to \infty} \displaystyle\sum_{n = 2}^\infty   t e^{\big(\lambda_1-\frac{j_n^2}{8(1-a)}\big)t} =0. 
 \]
 Taking into account this fact, we deduce from \eqref{eq7.th2} and \eqref{eq4.th2} that  there exists $T(u_0,\beta)>0$ such that the mapping $a\mapsto \partial_x u^a(1,t)$ is 
 increasing with respect to $a$ for all $t\ge T(u_0,\beta)$ and 
 
\begin{equation}\label{eq12.th2}
|\partial_a(\partial_x u^a(1,t))| \ge e^{-\lambda_1 t} \dfrac{L}{|j_1 J_0'(j_1)|}  \quad \forall\, t\ge T(u_0,\beta)
\end{equation}

 \medskip
 Finally let us obtain the stability estimate. We have for all $a_1,a_2\in [0,\beta]$
 
 \[
 \partial_x u^{a_2}(1,t) - \partial_x u^{a_1}(1,t) =  \displaystyle\sum_{n = 1}^\infty \frac{1}{j_n J'_0(j_n)}
 \Big[ f_n(t,a_2)U_n^0(a_2) - f_n(t,a_1)U_n^0(a_1)\Big].
 \]
 Therefore, in view of \eqref{eq12.th2} we conclude that  
 
 \[
| \partial_x u^{a_2}(1,t) - \partial_x u^{a_1}(1,t)| \ge e^{-\lambda_1 t} \dfrac{L}{|j_1 J_0'(j_1)|} |a_2-a_1|\quad \forall\, t\ge T(u_0,\beta),
 \]
 from where we obtain \eqref{eq.th} with $C = e^{\lambda_1 t} |j_1 J_0'(j_1)|/L $.

 Notice that if that $U_1^0(a)<- \delta \quad \forall\, a\in [0,\beta]$, then the argument is similar. 
 
 This ends the proof. \hfill $\blacksquare$
  \bigskip

 \section{Admissible initial data for one point measurement }\label{sec.examples}
 
 In this section we will present some examples of initial data for which we can have stability estimates with explicit  constants. Moreover, we will see that there exist  initial data for which Lipschitz stability estimates hold for an arbitrary small time.

 \medskip
 Let us first give the definition of an admissible initial data. 
 
\begin{defi}[Admissible initial data]\label{def41}
Let $I:=[t_0,t_1]\subseteq [0,+\infty)$, $K:=[\alpha,\beta]\subseteq [0,1)$. We will say that $u_0\in C^1([0,1])$ is an admissible initial value on $K\times I$ if there exists a positive constant $C$ such that 

\begin{equation}\label{eqstarBIS}
|\mu_{u_0}(a_1,t)- \mu_{u_0}(a_2,t) | \geq C |a_1 - a_2| \quad \forall\, a_1,a_2\in K,\quad   \forall\, t\in I,
\end{equation}
where $\mu_{u_0}$ is given by \eqref{eq.mapmu}. 
\end{defi}
 
 Notice that we are taking more regular initial data than ones for which the problem is well-posed.  
 
 \medskip
 In what follows, we will present some examples of admissible initial data.

 \bigskip
 \subsection{Case $u_0(x) = 1$}

We have the following result. 
 
\begin{figure}[h!]
 \centering 
 \includegraphics[width = 8cm]{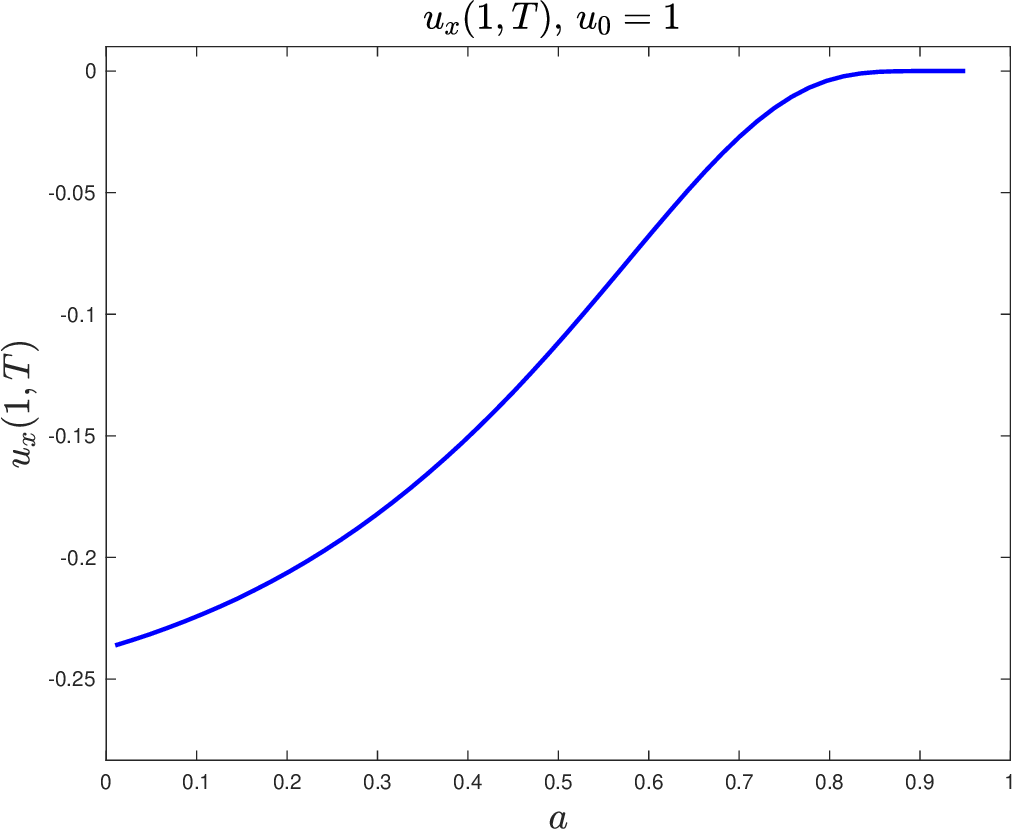}
 \caption{Case $u_0=1$: $\partial_x u^a(1,t)$ is strictly monotone with respect to $a$.}
 \label{fig.normder}
 \end{figure}

 \begin{prop}\label{prop.case1}
 Assume that $u_0=1$. We have that  $u_0$ is an admissible initial value in the sense of Definition \ref{def41}  for all $a_1,a_2\in [\alpha,\beta]\subseteq (0,1)$,  and all $t\in [t_0,t_1]\subseteq(0,+\infty)$. Moreover, the following Lipschitz stability estimate is fulfilled:
 
  \begin{equation}\label{eq.stabu01} 
 |a_2-a_1| \le   C |\partial_x u^{a_2}(1,t) - \partial_x u^{a_1}(1,t)| \quad \forall\, a_1,a_2\in [\alpha,\beta]\,\quad \forall\, t\in [t_0,t_1],
 \end{equation}
where  $u^{a_i}$, $i=1,2$ are solutions to \eqref{eq.pr} corresponding to $a_i$ and $C$ is given by 
\[
C = \displaystyle \dfrac{1}{t_0}e^{\big( \frac{j_1}{2}\big)^2 \frac{t_1}{1-\beta}}. 
\]
 \end{prop}
\medskip
\noindent
\textbf{Proof of Proposition~\ref{prop.case1}:}  In this case, from \eqref{eq.U0n} we have 
 
  \begin{equation*}
U^0_n(a)  =   \displaystyle \int_{0}^{j_n} u_0\Big(a+(1-a)\frac{s^2}{j_n^2} \Big)s J_0 (s) \, ds
= \displaystyle \int_{0}^{j_n} s J_0 (s) \, dx = -j_n J_0'(j_n). 
 \end{equation*}
 Therefore, by \eqref{eq.normderuxa} we have 
 \begin{equation}\label{eq.normderuxaB}
\partial_x u^a(1,t)  =  -\displaystyle\sum_{n = 1}^\infty  f_n(t,a) =   
 -\displaystyle\sum_{n = 1}^\infty e^{-\big( \frac{j_n}{2}\big)^2\frac{t}{1-a}},
 \end{equation} 
 which is strictly monotone with respect to $a$ (see Figure~\ref{fig.normder}, where we have taken $t = T = 1$). Therefore, in particular, in this case, we have the uniqueness for our inverse problem. 
 
 \medskip

 \medskip
 \noindent 
 For all $0<a_1\le a_2 <1$, we also obtain 
 
 \[
 \begin{split}
 \partial_x u^{a_2}(1,t) - \partial_x u^{a_1}(1,t)  & = \displaystyle\sum_{n = 1}^\infty \big [f_n(t,a_1)- f_n(t,a_2) \big] 
 \\[3mm]
& \ge   f_1(t,a_1) -  f_1(t,a_2) =  f_1(t,a_2) \Big( \dfrac{f_1(t,a_1)}{f_1(t,a_2)}-1\Big) 
 \\[3mm]
& \ge
 f_1(t,a_2) \Big(e^{-\big( \frac{j_1}{2}\big)^2 t \frac{a_1 - a_2}{(1-a_1)(1-a_2)}}-1\Big)\ge f_1(t,a_2) \big( \frac{j_1}{2}\big)^2 t (a_2-a_1).
 \end{split}
 \]
Therefore, we have the following Lipschitz stability estimate: 
 
 \begin{equation*}
 |a_2-a_1| \le   \Big( \frac{2}{j_1^2}\Big)^2  \displaystyle \dfrac{1}{t f_1(t,a_2)}  |\partial_x u^{a_2}(1,t) - \partial_x u^{a_1}(1,t)|
 \le  
 \dfrac{1}{t_0} e^{\big( \frac{j_1}{2}\big)^2 \frac{t_1}{1-\beta}}  |\partial_x u^{a_2}(1,t) - \partial_x u^{a_1}(1,t)|.
 \end{equation*}
 \noindent
 This ends the proof. \hfill $\blacksquare$

 \bigskip
 \subsection{Case $u_0(x) = 1-x$}

We have the following Lipschitz stability result. 
 
 \begin{prop}\label{prop.case2}
 Assume that $u_0=1-x$. We have that  $u_0$ is an admissible initial value in the sense of Definition \ref{def41}  for all $a_1,a_2\in [\alpha,\beta]\subseteq (0,1)$,  and all $t\in [t_0,t_1]\subseteq[0,+\infty)$. Moreover, the following Lipschitz stability estimate is fulfilled:
 
  \begin{equation}\label{eq.stabu01Case2} 
 |a_2-a_1| \le   C |\partial_x u^{a_2}(1,t) - \partial_x u^{a_1}(1,t)| \quad \forall\, a_1,a_2\in [\alpha,\beta]\,\quad \forall\, t\in [t_0,t_1],
 \end{equation}
where  $u^{a_i}$, $i=1,2$ are solutions to \eqref{eq.pr} corresponding to $a_i$ and $C$ is given by 
\[
C = \displaystyle \dfrac{\pi^2}{16}e^{\big( \frac{j_1}{2}\big)^2 \frac{t_1}{1-\beta}}. 
\]

  \end{prop}
\medskip
\noindent
\textbf{Proof of Proposition~\ref{prop.case2}:}  Taking into account properties from Lemma~\ref{lemma.Bessel}, we have 
 
 \begin{equation}\label{eq.Uaex2}
 \begin{split}
U^0_n(a)  & =  \displaystyle \int_{0}^{j_n} (1-a)\Big(1-\dfrac{s^2}{j_n^2} \Big)s J_0 (s) \, ds 
= 
\displaystyle (1-a)\Bigg[ \int_{0}^{j_n} s J_0 (s) \, ds - \int_{0}^{j_n} \dfrac{s^3}{j_n^2} J_0 (s) \, ds\Bigg]
\\[3mm]
& = \displaystyle (1-a)\Bigg[ -j_nJ'_0(j_n) 
- 
\frac{1}{j_n^2}\Big\{ s^2 s J_1(s) \Big|_{0}^{j_n} + 2 \int_{0}^{j_n} s^2 J_1 (s) \, ds \Big\}\Bigg]
\\[3mm]
& = \displaystyle (1-a)\Bigg[ -j_nJ'_0(j_n) 
- 
\dfrac{1}{j_n^2} j_n^3 J_1(j_n)+ \dfrac{2}{j_n^2}  s^2 J_2(s)\Big|_{0}^{j_n} \Bigg]
\\[3mm]
& = \displaystyle (1-a)\Bigg[ -j_nJ'_0(j_n) -  j_n J_1(j_n)+ 2J_2(j_n) \Bigg] = 2(1-a) J_2(j_n)
\\[3mm]
& = 2\displaystyle (1-a) \Big(\dfrac{2}{j_n} J_1(j_n)-J_{0}(j_n)\Big) = -\dfrac{4(1-a)}{j_n} J_0'(j_n). 
\end{split} 
 \end{equation}
 
 \noindent 
Here, we have used the fact that $J_0(j_n) =0$ and property \textit{b)} form Lemma~\ref{lemma.Bessel}, which implies that  $J_2(j_n) = \dfrac{2}{j_n} J_1(j_n)-J_{0}(j_n)$.  
 
\begin{figure}[h!]
 \centering 
 \includegraphics[width = 8cm]{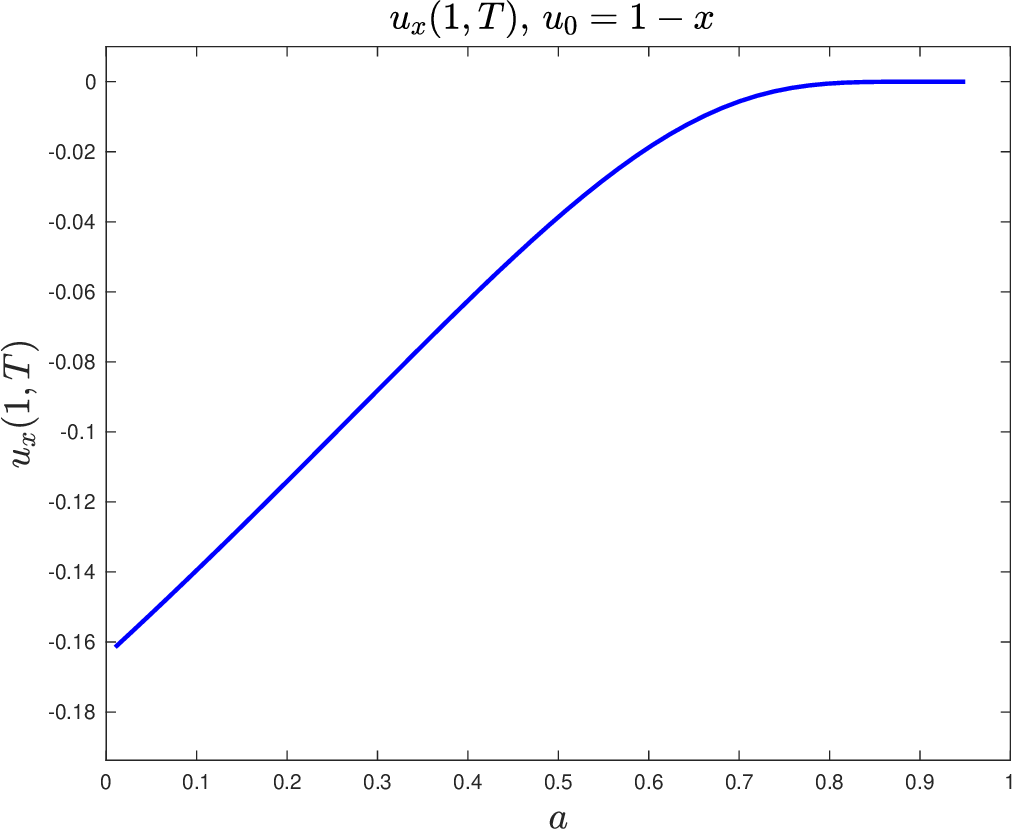}
 \caption{Case $u_0=1-x$: $\partial_x u^a(1,t)$ is strictly monotone with respect to $a$.}
 \label{fig.normderux}
 \end{figure}

Therefore, using the above expression in \eqref{eq.normderuxa}, we obtain 
 \begin{equation*}
\partial_x u^a(1,t)  = -4(1-a)  \displaystyle\sum_{n = 1}^\infty \dfrac{f_n(t,a)}{j_n^2}, 
 \end{equation*} 
 which is again monotone with respect to $a$ (see Figure~\ref{fig.normderux}, where we have taken $t = T = 1$), since 
 
 \begin{equation}\label{normdercase2}
 \partial_a( (1-a) f_n(t,a) )= -f_n(t,a) + (1-a) \partial_a f_n(t,a) <0.  
 \end{equation}
   \medskip

\medskip
\noindent
Moreover, for all $0<a_1\le a_2 <1$, we have that
 
 \[
 \begin{split}
 \partial_x u^{a_2}(1,t) - \partial_x u^{a_1}(1,t)  & 
 = \displaystyle \frac{4}{j_1^2}\Big [ (1-a_1)\sum_{n = 1}^\infty f_n(t,a_1)- (1-a_2)\sum_{n = 1}^\infty f_n(t,a_2) \Big] 
 \\[3mm]
& \ge \frac{4}{j_1^2}(  (1-a_1) f_1(t,a_1) - (1-a_2) f_1(t,a_2) ),  
 \end{split}
 \]
from which, using \eqref{normdercase2} we obtain, for some $\bar{a} \in (a_1,a_2)$,

 \[
 \begin{split}
 |\partial_x u^{a_2}(1,t) - \partial_x u^{a_1}(1,t) | & 
 \ge \displaystyle  \frac{4}{j_1^2} \Big |\partial_a( (1-a) f_1(t,a) \Big|_{a=\bar{a}}) \Big| |a_2-a_1|
 \\[3mm]
& \ge\frac{4}{j_1^2} \min_{a\in[\alpha,\beta], t\in[t_0,t_1]} \Big[ \Big( 1+\big(\frac{j_1}{2}\big)^2 \dfrac{t}{1-a}\Big) f_1(t,a)\Big]|a_2-a_1|.
 \\[3mm]
& \ge
\frac{4}{j_1^2}  e^{-\big( \frac{j_1}{2}\big)^2 \frac{t_1}{1-\beta}}    |a_2-a_1|
\ge 
\frac{16}{\pi^2}  e^{-\big( \frac{j_1}{2}\big)^2 \frac{t_1}{1-\beta}}    |a_2-a_1|, 
 \end{split}
 \]
 since $1/j_1^2 \ge 4/\pi^2$. This ends the proof. \hfill $\blacksquare$

 \bigskip\bigskip
 \subsection{Case $u_0(x) = x(1-x)$}
   
The following holds. 

 \begin{prop}\label{prop.case3}
 Assume that $u_0=x(1-x)$. We have that  $u_0$ is an admissible initial value in the sense of Definition \ref{def41}  for all $a_1,a_2\in [\alpha,\beta]\subseteq (0,1)$,  and all $t\in [t_0,t_1]\subseteq[0,+\infty)$ with $t_0$ large enougnt. Moreover, the following Lipschitz stability estimate is fulfilled:
 
  \begin{equation}\label{eq.stabu01Case3} 
 |a_2-a_1| \le   C |\partial_x u^{a_2}(1,t) - \partial_x u^{a_1}(1,t)| \quad \forall\, a_1,a_2\in [\alpha,\beta]\,\quad \forall\, t\in [t_0,t_1],
 \end{equation}
where  $u^{a_i}$, $i=1,2$ are solutions to \eqref{eq.pr} corresponding to $a_i$ and a positive constant $C$. 
\end{prop}

\medskip
\noindent
\textbf{Proof of Proposition~\ref{prop.case3}:}

   \begin{equation}\label{eq.UaEj3}
 \begin{split}
U^0_n(a)  & =  \displaystyle \int_{0}^{j_n} \Big(a+(1-a) \dfrac{s^2}{j_n^2} \Big) (1-a)\Big(1-\dfrac{s^2}{j_n^2} \Big)s J_0 (s) \, ds \\[3mm]
& = 
a(1-a) \displaystyle \int_{0}^{j_n} \Big(1-\dfrac{s^2}{j_n^2} \Big)s J_0 (s) \, ds 
+
(1-a)^2 \displaystyle \int_{0}^{j_n} \dfrac{s^2}{j_n^2} \Big(1-\dfrac{s^2}{j_n^2} \Big)s J_0 (s) \, ds 
\\[3mm]
&
:=
V_n(a) + W_n(a). 
\end{split} 
 \end{equation}
Using \eqref{eq.Uaex2}, we obtain 
 
 \begin{equation}\label{eq.VaEj3}
 V_n(a) = a(1-a) \displaystyle \int_{0}^{j_n} \Big(1-\dfrac{s^2}{j_n^2} \Big)s J_0 (s) \, ds = -\dfrac{4a(1-a)}{j_n} J_0'(j_n). 
 \end{equation}

\noindent
On the other hand, we can wtrire that 

 \begin{equation}\label{eq.WaEj3}
 \begin{split}
 W_n(a) & = (1-a)^2 \displaystyle \int_{0}^{j_n} \dfrac{s^2}{j_n^2} \Big(1-\dfrac{s^2}{j_n^2} \Big)s J_0 (s) \, ds 
 \\[3mm]
 & =
 (1-a)^2 \Big[  \dfrac{1}{j_n^2} \displaystyle \int_{0}^{j_n} s^3 J_0 (s) \, ds 
 - 
  \dfrac{1}{j_n^4} \displaystyle \int_{0}^{j_n} s^5 J_0 (s) \, ds. 
 \Big]  
 \end{split}
 \end{equation}
Now,  arguing as in \eqref{eq.Uaex2},  we deduce

 \begin{equation}\label{eq.WaEj3B}
\begin{split}
 \dfrac{1}{j_n^2} \displaystyle \int_{0}^{j_n} s^3 J_0 (s) \, ds  
 & =  \dfrac{1}{j_n^2} \Big( j_n^3 J_1(j_n) - 2 j_n J_2(j_n)\Big) = j_n J_1(j_n) - 2 J_2(j_n)
 \\[3mm] 
& = -j_n J'(j_n) - \dfrac{4}{j_n} J_1(j_n) 
  = \Big(\dfrac{4}{j_n} - j_n \Big) J_0'(j_n) 
 \end{split} 
 \end{equation}
 Also, we have 
 
 \begin{equation}\label{eq.s5J}
\begin{split}
 \dfrac{1}{j_n^4} \displaystyle \int_{0}^{j_n} s^5 J_0 (s) \, ds  
 & =  \dfrac{1}{j_n^4} \Big[ s^4 s J_1(j_n)\Big|_{0}^{j_n} - 4 \displaystyle \int_{0}^{j_n} s^4 J_1 (s) \, ds \Big] 
 \\[3mm] 
&
= 
 j_n J_1(j_n) - \dfrac{1}{j_n^4} \Big[ s^2 s^2 J_2(j_n)\Big|_{0}^{j_n} - 2 \displaystyle \int_{0}^{j_n} s^3 J_2 (s) \, ds \Big] 
 \\[3mm] 
&  = 
 j_n J_1(j_n) -  4 J_2(j_n) + \dfrac{8}{j_n} J_3(j_n). 
 \end{split} 
\end{equation}

 Moreover, thanks to the property \textit{b)} from Lemma~\ref{lemma.Bessel}, we have 
 
 \begin{equation}\label{eq.J2}
 J_2(j_n) = \dfrac{2}{j_n} J_1(j_n) - J_0(j_n) =  \dfrac{2}{j_n} J_1(j_n)
\end{equation}
and 

 \begin{equation}\label{eq.J3}
 J_3(j_n) = \dfrac{4}{j_n} J_2(j_n) - J_1(j_n) =  \dfrac{8}{j_n^2} J_1(j_n)- J_1(j_n) = \Big(  \dfrac{8}{j_n^2} -1 \Big) J_1(j_n) 
\end{equation}
Therefore, taking onto account \eqref{eq.J2} and \eqref{eq.J3}  in  \eqref{eq.s5J}, we deduce 

 \begin{equation}\label{eq.s5Jfin}
\begin{split}
 \dfrac{1}{j_n^4} \displaystyle \int_{0}^{j_n} s^5 J_0 (s) \, ds  
 & =
 j_n J_1(j_n) -  \dfrac{8}{j_n} J_1(j_n) + \dfrac{8}{j_n} \Big(  \dfrac{8}{j_n^2} -1 \Big) J_1(j_n) 
 \\[3mm]
& =  j_n\Big( 1- \dfrac{8}{j_n^2} \Big)^2 J_1(j_n) = -j_n\Big( 1- \dfrac{8}{j_n^2} \Big)^2 J_0'(j_n).
 \end{split} 
\end{equation}

By combining \eqref{eq.UaEj3}, \eqref{eq.VaEj3}, \eqref{eq.WaEj3}, \eqref{eq.WaEj3B} and \eqref{eq.s5Jfin} we finally obtain

\[
\begin{split}
U^0_n(a)  & = -\dfrac{4a(1-a)}{j_n} J_0'(j_n) + (1-a)^2 \Big(\dfrac{4}{j_n} - j_n \Big) J_0'(j_n) +(1-a)^2  j_n\Big( 1- \dfrac{8}{j_n^2} \Big)^2 J_0'(j_n) 
\\[3mm]
& = J_0'(j_n) \Big[  -\dfrac{4a(1-a)}{j_n}  + (1-a)^2  j_n\Big\{ \dfrac{4}{j_n^2} - 1 + \Big( 1- \dfrac{8}{j_n^2} \Big)^2\Big\}
\Big] := J_0'(j_n) K(a,j_n), 
\end{split}
\]
where  $K(a,j_n)$ is given by

\begin{equation}\label{eq.K}
K(a,j_n) =  -\dfrac{4a(1-a)}{j_n}  + (1-a)^2  \dfrac{4}{j_n}  \Big(\dfrac{16}{j_n^2}-3 \Big).
\end{equation}

\noindent
Therefore, from \eqref{eq.normderuxa} we obtain 
 \begin{equation*}
\partial_x u^a(1,t)  = \displaystyle\sum_{n = 1}^\infty \dfrac{f_n(t,a)}{j_n}K(a,j_n). 
 \end{equation*} 
It is not difficult to see that $\partial_x u^a(1,t)$ is not strictly monotone (see Figure~\ref{fig.normderux1}). Indeed, using \eqref{eq.K}  and  that 

\[
\partial_a K(a,j_n) = - \dfrac{4}{j_n}(1-2a) + 2(a-1)\dfrac{4}{j_n}  \Big( \dfrac{16}{j^2_n} -3 \Big),
\]
we have 

 \begin{equation}\label{eq.2der}
\begin{split} 
\partial^2_{ax} u^a(1,t)  &  = \displaystyle\sum_{n = 1}^\infty \dfrac{1}{j_n}\Big[    \partial_a f(t,a) K(a,j_n) + f_n(t,a) \partial_a K(a,j_n) \Big]  
\\[2mm]
&  = \displaystyle\sum_{n = 1}^\infty \dfrac{1}{j_n}\Big[    -\Big(\dfrac{j_n}{2} \Big)^2 \dfrac{t}{(1-a)^2}f_n(t,a) K(a,j_n) + f_n(t,a) \partial_a K(a,j_n) \Big]  
\\[2mm]
& = \displaystyle\sum_{n = 1}^\infty \dfrac{f_n(t,a) }{j^2_n}
\Bigg[
\dfrac{j_n^2 a}{1-a}t - t(16-3j_n^2) - 4(1-2a) + 8 (a-1)\Big( \dfrac{16}{j^2_n} -3 \Big)
 \Bigg]. 
\end{split}
 \end{equation} 

This shows that $\partial^2_{ax} u^a(1,t) >0$ for $t$ sufficiently large, as guaranteed by Theorem~\ref{th.stabillity}. Therefore, the map $a\mapsto \partial_x u^a(1,t)$ is strictly increasing. 

On the other hand, we claim that   $a\mapsto \partial_x u^a(1,t)$ fails to be monotone for $t>0$ sufficiently small. Indeed, 

 \begin{equation*}
\begin{split} 
\lim_{t\to 0}\partial^2_{ax} u^a(1,t) \Big|_{a= 0} 
&  = \lim_{t\to 0} -  \displaystyle\sum_{n = 1}^\infty \dfrac{f_n(t,0) }{j^2_n}
\Big[ t(16-3j_n^2) + 4 + 8 \Big( \dfrac{16}{j^2_n} -3 \Big) \Big]  
 \\[3mm]
& = 
-  \displaystyle\sum_{n = 1}^\infty \dfrac{1}{j^2_n}
\Big[ 4 + 8 \Big( \dfrac{16}{j^2_n} -3 \Big) \Big]  < 0, 
\end{split}
 \end{equation*} 

\noindent 
because the above function series converges uniformly on $[0,\infty)$. Consequently,  we obtain $\partial^2_{ax} u^a(1,t)\Big |_{a= 0}<0$ for $t$ sufficiently small. Moreover, the term in the square brackets in the expression  \eqref{eq.2der} is positive for  $a$ sufficiently close to 1, say $a>\alpha$. Therefore, $\partial^2_{ax} u^a(1,t)>0$ for all such values of $a$. 
   
  \medskip  
 \begin{figure}[h!]
 \centering 
 \includegraphics[width = 7cm]{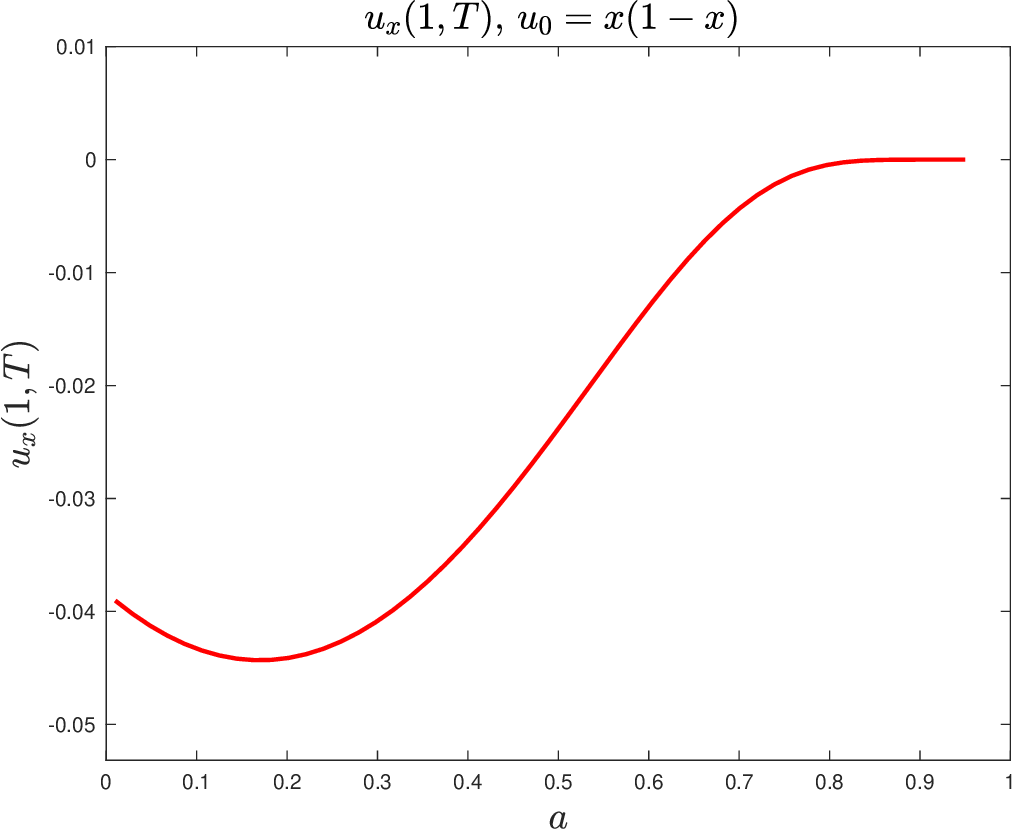}
 \includegraphics[width = 6.8cm]{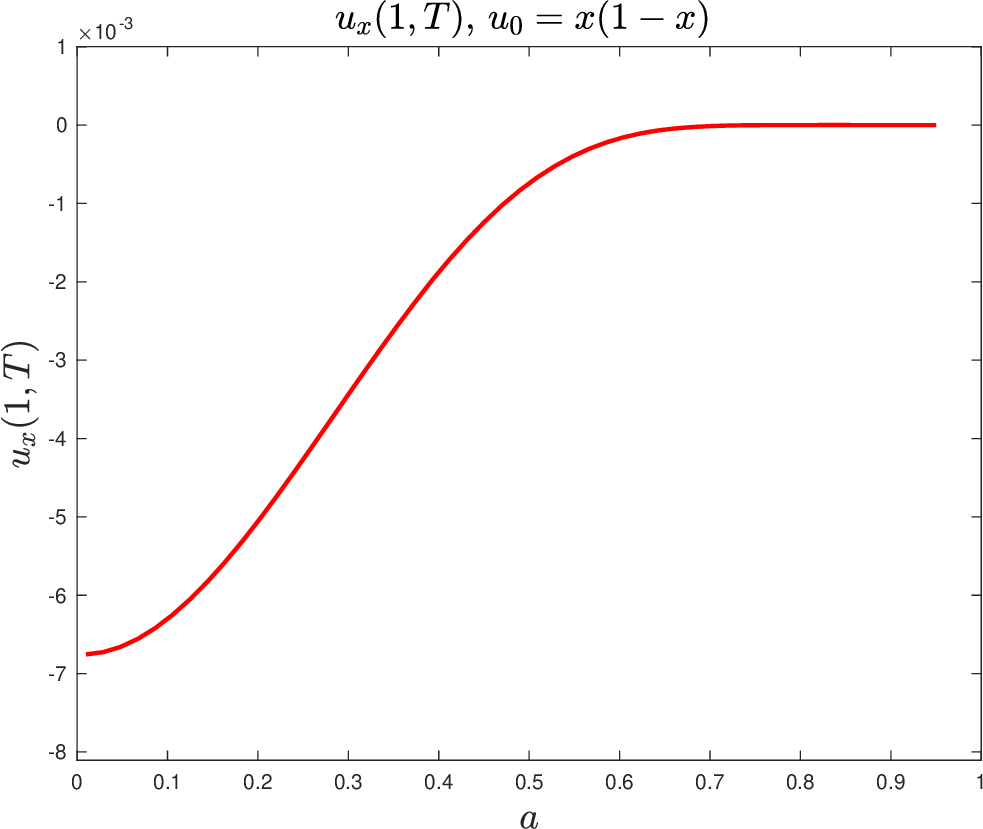}
 \caption{Case $u_0=x(1-x)$: $\partial_x u^a(1,t)$ is not monotone with respect to $a$ (left, $T=1$), $\partial_x u^a(1,t)$ is monotone for $T$ large enough (right, $T = 2.2$).}
 \label{fig.normderux1}
 \end{figure}
 
 Finally, the estimate \eqref{eq.stabu01Case3} can be deduced using similar arguments as in  previous examples.  This ends the proof. 
\hfill $\blacksquare$ 

\bigskip

The above reasoning explains  the behavior observed in the Figure~\ref{fig.normderux1} below. 
 \bigskip
 \eject
 \subsection{Case $u_0(x) = x$}

In this case, the following is found. 
 \begin{prop}\label{prop.case4}
 Assume that $u_0=x$. We have that  $u_0$ is an admissible initial value in the sense of Definition \ref{def41}  for all $a_1,a_2\in [\alpha,\beta]\subseteq (0,1)$,  and all $t\in [t_0,t_1]\subseteq[0,+\infty)$ with $t_0$ large enougnt. Moreover, the following Lipschitz stability estimate is fulfilled:
 
  \begin{equation*}
 |a_2-a_1| \le   C |\partial_x u^{a_2}(1,t) - \partial_x u^{a_1}(1,t)| \quad \forall\, a_1,a_2\in [\alpha,\beta]\,\quad \forall\, t\in [t_0,t_1],
 \end{equation*}
where  $u^{a_i}$, $i=1,2$ are solutions to \eqref{eq.pr} corresponding to $a_i$ and $C$ is given by 
\[
C = \displaystyle j_1^2 e^{\big( \frac{j_1}{2}\big)^2 \frac{t_1}{1-\beta}}.
\]
\end{prop}
\medskip
\noindent
\textbf{Proof of Proposition~\ref{prop.case4}:} It is not difficult to see that we can proceed as in the proof of Theorem~	\ref{th.stabillity}. Let us just notice that for $a\in [\alpha,\beta]\subseteq [0,1)$ and $t$ sufficiently large we have the representation of $\partial_x u^a(1,t)$ given in Figure~\ref{fig.normderux1} (right). 
%

 \begin{figure}[h!]
 \centering 
 \includegraphics[width =7cm]{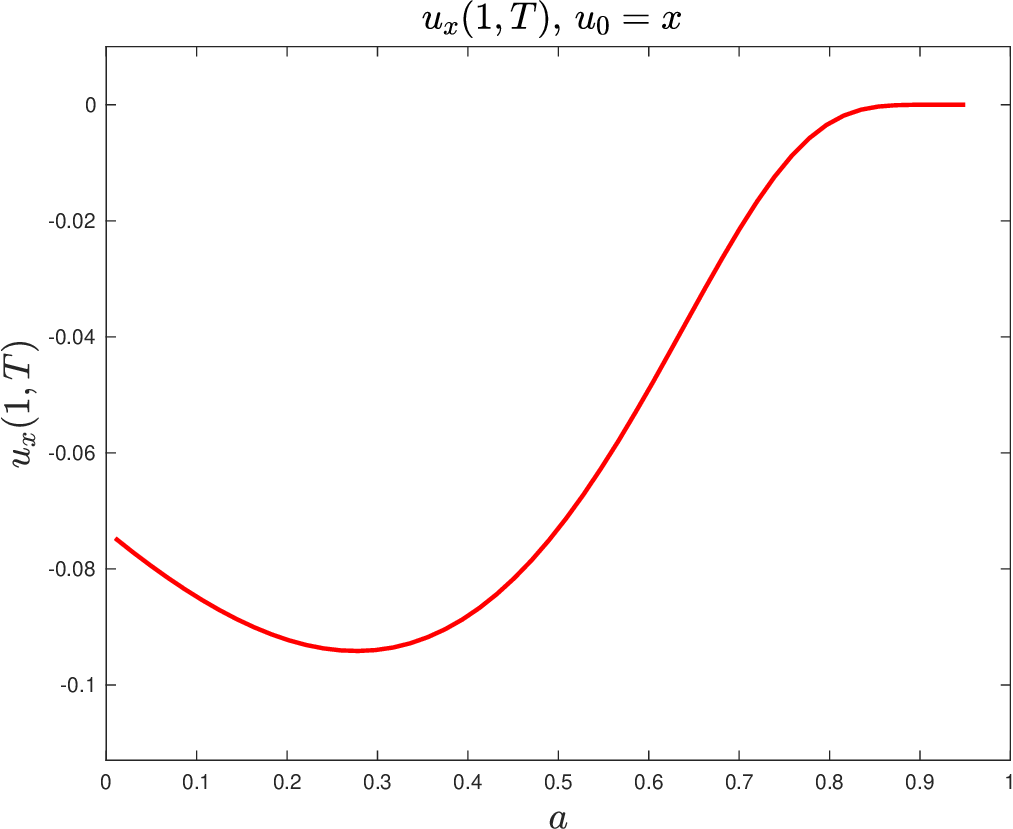}
  \includegraphics[width = 7cm]{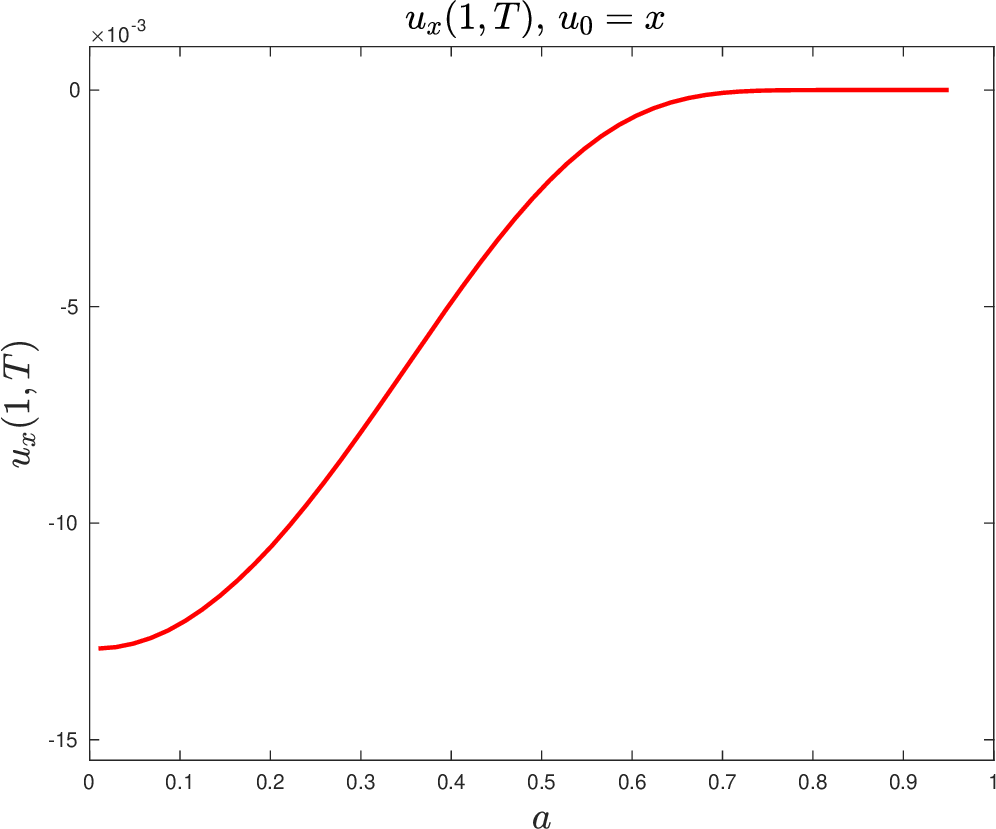}
 \caption{Case $u_0=x$: $\partial_x u^a(1,t)$ is not monotone in $a$ (left, $T=1$), $\partial_x u^a(1,t)$ is monotone for $T$ large enough (right, $T = 2.2$).}
 \label{fig.normderux1}
 \end{figure}

\hfill $\blacksquare$

\bigskip
\section{Uniqueness results for ``distributed'' measurements}\label{sec:uniqueness}

In this section we will present general uniqueness results for  \eqref{eq.pr} based on the explicit representation given in Theorem~\ref{prop.normder}.  In this case, in contrast to the previous sections where we consider point wise measurements, we will need measurements distributed over a  time interval.

 \begin{theorem}\label{theorem2}
Let $0<a_1,a_2 <1$ and $0<t_1<t_2$.  Let   $u^{a_1}$ and $u^{a_2}$ be the solutions to \eqref{eq.pr} corresponding to initial values $u_0$ and $\widetilde{u}_0$, respectively. Assume that $U^0_1(a_1)$ and $\widetilde{U}^0_1(a_2)$  (given by~\eqref{eq.U0n} for $u_0$ and $\widetilde{u}_0$, respectively) satisfy

\begin{equation}\label{th2.eq1}
U_1^0(a_1) \neq 0 \quad \text{and} \quad \widetilde{U}_1^0(a_2) \neq 0. 
\end{equation}
\noindent 
Then $\partial_x u^{a_1}(1,t) = \partial_x u^{a_2}(1,t)$ for $t_1<t<t_2$ implies that $a_1=a_2$ and $u_0 = \widetilde{u}_0$. 

 \end{theorem}
\medskip 
\textbf{Proof of Theorem~\ref{theorem2}:} Note that the function $t\mapsto \partial_x u^a(1,t)$ is analytic for all $t>0$.  Let us set 

\begin{equation}\label{lambdamu}
\lambda_n:= \dfrac{j^2_n}{4(1-a_1)}, \quad \mu_n:= \dfrac{j^2_n}{4(1-a_2)}, \quad n\in \mathds{N}. 
\end{equation}

\noindent
Notice that $\lambda_1<\lambda_2<\cdots$ and $\mu_1<\mu_2<\cdots$.  Then the time-analyticity and \eqref{eq.normderuxa} yield 

\begin{equation}\label{eq2.uniq2t2}
\displaystyle\sum_{n=1}^\infty \dfrac{U_n^0(a_1)}{j_n J'_0(j_n)} e^{-\lambda_n t} 
=
\displaystyle\sum_{n=1}^\infty \dfrac{\widetilde{U}_n^0(a_2)}{j_n J'_0(j_n)} e^{-\mu_n t} , \quad t>t_1 \quad (t>0). 
\end{equation}

\noindent
Hence, 

\[
\dfrac{U_1^0(a_1)}{j_1 J'_0(j_1)} e^{-\lambda_1 t} + \displaystyle\sum_{n=2}^\infty \dfrac{U_n^0(a_1)}{j_n J'_0(j_n)} e^{-\lambda_n t} 
=
\dfrac{\widetilde{U}_1^0(a_2)}{j_1 J'_0(j_1)} e^{-\mu_1 t} + \displaystyle\sum_{n=2}^\infty \dfrac{\widetilde{U}_n^0(a_2)}{j_n J'_0(j_n)} e^{-\mu_n t} , \quad t>t_1. 
\]

\medskip
Assume that $a_1\neq a_2$. Without  loss of generality, we can assume that $a_1<a_2$, that is $\frac{1}{1-a_1}<\frac{1}{1-a_2}$. Then $\lambda_n<\mu_n$ for $n\in \mathds{N}$. Therefore, we have 

\begin{equation}\label{eq.uniq2t2}
\begin{split}
\dfrac{U_1^0(a_1)}{j_1 J'_0(j_1)} & + \displaystyle\sum_{n=2}^\infty \dfrac{U_n^0(a_1)}{j_n J'_0(j_n)} e^{-(\lambda_n-\lambda_1) t} \\[3mm]
 & =
\dfrac{\widetilde{U}_1^0(a_2)}{j_1 J'_0(j_1)}e^{-(\mu_1-\lambda_1) t}+ \displaystyle\sum_{n=2}^\infty \dfrac{\widetilde{U}_n^0(a_2)}{j_n J'_0(j_n)} e^{-(\mu_n-\lambda_1) t}, \quad t>t_1. 
\end{split}
\end{equation}

\noindent
Since $\mu_n-\lambda_n>0$, we let $t\to\infty$ in the previous equality to obtain

\[
\dfrac{U_1^0(a_1)}{j_1 J'_0(j_1)} = 0,
\] 
in contrast with \eqref{th2.eq1}.  Thus $a_1=a_2$ follows.

Moreover, the above argument also yields $\lambda_n=\mu_n$ for all $n\in \mathds{N}$.

\medskip
Let us now see that $u_0 = \widetilde{u}_0$. Indeed, since $\lambda_n=\mu_n$ for all $n\in \mathds{N}$ and $a_1=a_2\equiv a$, form \eqref{eq2.uniq2t2} we have 

\begin{equation}\label{eq.u0}
\displaystyle\sum_{n=1}^\infty \dfrac{e^{-\lambda_n t}}{j_n J'_0(j_n)} \big( U_n^0(a)-\widetilde{U}_n^0(a)\big)= 0 , \quad t>t_1 \quad (t>0). 
\end{equation}

Let us set 

\[
n_0 = \inf\{ n\ge 1:  U_{n_0}^0(a)\neq \widetilde{U}_{n_0}^0(a)  \}.
\]

We are going to show that this is an empty set or, equivalently, $n_0=\infty$. Suppose $n_0<\infty$ and multiply equality \eqref{eq.u0} by $e^{\lambda_{n_0}t}$   to obtain 

\begin{equation*}
\displaystyle\sum_{n=1}^\infty \dfrac{e^{-\lambda_n t}}{j_n J'_0(j_n)} \big( U_n^0(a)-\widetilde{U}_n^0(a)\big)
e^{\lambda_{n_0} t}= 0, \quad t>t_1
\end{equation*}

\noindent
and then

\begin{equation*}
U_{n_0}^0(a)-\widetilde{U}_{n_0}^0(a)  +\displaystyle\sum_{n=n_0+1}^\infty \dfrac{e^{-(\lambda_n-\lambda_{n_0})t}}{j_n J'_0(j_n)} \big( U_n^0(a)-\widetilde{U}_n^0(a)\big)= 0, \quad t>t_1.
\end{equation*}
\noindent
We let $t\to +\infty$ to deduce from the previous equality that $U_{n_0}^0(a)=\widetilde{U}_{n_0}^0(a)$ in contrast with with the definition of $n_0$. Therefore, $n_0$ is empty and $U_{n}^0(a)=\widetilde{U}_{n}^0(a)$ for all $n$. 

From \eqref{eq.solua} we conclude that $u_0 =\widetilde{u}_0$ by the coincidence of all the Fourier coefficients. This ends the proof.

\hfill$\blacksquare$

\medskip

We have the following consequence of this theorem.

 \begin{coro}\label{cor.theorem2}
Let   $u^{a_1}$ and $u^{a_2}$ be the solutions to \eqref{eq.pr} corresponding to initial values $u_0$ and $\widetilde{u}_0$, respectively. We assume 
\begin{equation*}
u_0\ge 0, \, u_0\neq 0 \text{ in } (a_1,1), \qquad \widetilde{u}_0\ge 0, \, \widetilde{u}_0\neq 0 \text{ in } (a_2,1)
\end{equation*}

Then $\partial_x u^{a_1}(1,t) = \partial_x u^{a_2}(1,t)$ for $t_1<t<t_2$ implies $a_1=a_2$ and $u_0=\widetilde{u}_0$. 

 \end{coro}
\medskip
\noindent
\textbf{Proof of Corollary~\ref{cor.theorem2}:} Since $J_0(s)>0$ for $0<s<j_1$, and $u_0\ge 0$, $u_0\neq 0$ in $(a_1,1)$ and 
$ \widetilde{u}_0\ge 0$, $\widetilde{u}_0\neq 0$ in $(a_2,1)$, we see $U_1^0(a_1)>0$ and $U_1^0(a_2)>0$. Thus, by Theorem ~\ref{theorem2} one completes the proof. 

\hfill$\blacksquare$

\begin{remark}
Notice that we have similar results for the solution $w^{a}$ to \eqref{eq.cp} corresponding to initial values $w_0$ with $w_0 = u_0$ in $(a,1)$ and $w_0 = v_0$ in $(0,a)$.  We can proceed to the consideration of the determination of $a$ and the initial values. Indeed if  $w^{a_1}$ and $w^{a_2}$ are the solutions to \eqref{eq.cp} corresponding to $a_1$ and $a_2$ and the initial values $(u_0,v_0)$ and 
$(\widetilde{u}_0,\widetilde{v}_0)$, respectively, then we can conclude that $a_1=a_2$ and $u_0 = \widetilde{u}_0$ and  $v_0 = \widetilde{v}_0$.
\end{remark}

Now, let us present some characterization for the degeneracy points yielding the same observation data.

\begin{theorem}\label{th.uniq2}
We assume
$$
\ppp_xu^{a_1}(1,t) = \ppp_xu^{a_2}(1,t) \quad \mbox{for $t_1 < t < t_2$},
$$
and $u^{a_1} \not\equiv 0$ or $u^{a_2} \not\equiv 0$ in $(0,1) \times (0,T)$.
Then there exist $m_1,m_2 \in \mathds{N}$ such that 

\begin{equation}\label{5.1}
\frac{1-a_1}{j_{m_1}^2} = \frac{1-a_2}{j_{m_2}^2}.
\end{equation}

\end{theorem}

\begin{remark}
We can verify that, if $m_1=m_2$, then $a_1=a_2$, that is, the uniqueness in determining the degeneracy point holds.

Since $(u^{a_1}, u^{a_2}) \equiv (0,0)$ is a trivial case,
we can interpret that the non-uniqueness in determining $a_1$, $a_2$ essentially occurs only in the case of 
\eqref{5.1} with $m_1 \ne m_2$.   Since $j_n$, $n\in \mathds{N}$ is known to be 
transcendental numbers (see \cite{Siegel}), this means that uniqueness occurs only if
$$
\mbox{$\sqrt{\frac{1-a_1}{1-a_2}} \ne 1$ is 
given by a ratio of two transcendental numbers}.
$$

\noindent
Thus, non-uniqueness happens only in special locations of degeneracy.
\end{remark}

\medskip
For the proof of Theorem~\ref{th.uniq2}, we we need use the following result:

\begin{lemma}\label{lemma5.1}
Let $\{b_n\}$ be a real sequence such that 

\begin{equation}\label{eq2Lemma5.1}
\sum_{n=1}^\infty \frac{b_n^2}{n^{2\gamma}} < \infty
\end{equation}
and let
$$
0 < \theta_1 < \theta_2 < .... \quad \mbox{and}\quad
c_0n^2 + o(n^2) \le \theta_n \quad \mbox{as $n\to\infty$}
$$
for some constant $\gamma > 0$ and $c_0 > 0$. Let $t_1>0$ be arbitrarily chosen. 
Then 
\begin{equation}\label{eqLemma5.1}
\sum_{n=1}^{\infty} b_ne^{-\theta_nt} = 0 \quad \mbox{for almost all 
$t>t_1$},
\end{equation}
implies $b_n=0$ for all $n\in \mathds{N}$.

\medskip

\end{lemma}
\noindent
{\bf Proof of Lemma~\ref{lemma5.1}:} Setting $r_n:= \theta_n - \theta_1 > 0$ and multiplying \eqref{eqLemma5.1} by $e^{\theta_1 t}$ we obtain
\[
b_1 + \sum_{n=2}^{\infty} b_ne^{-r_nt} = 0, \quad t > t_1.
\]
Then $\displaystyle\lim_{t\to \infty} b_ne^{-r_nt} = 0$ for all $n\in \mathds{N}$ and taking into account \eqref{eq2Lemma5.1} we have
\begin{equation}
\sum_{n=2}^{\infty} \vert b_ne^{-r_nt} \vert 
= \sum_{n=2}^{\infty} \left\vert \frac{b_n}{n^{\gamma}} n^{\gamma}e^{-r_nt}
\right\vert 
\le \left( \sum_{n=2}^{\infty}\left\vert \frac{b_n^2}{n^{2\gamma}}\right\vert 
\right)^{\frac{1}{2}} 
\left(\sum_{n=2}^{\infty} n^{2\gamma}e^{-2r_nt}\right)^{\frac{1}{2}}< \infty,
\end{equation}
because 
$$
n^{2\gamma}e^{-2r_nt} \le n^{2\gamma}e^{-2r_nt_2}
\le C_2n^{2\gamma}e^{-2c_0n^2t_2(1+o(1))}
\le C_3n^{2\gamma}e^{-C_3n^2}
$$
as $n\to \infty$ with some positives constants $C_2$ and $C_3$.

Therefore we see
$$
\lim_{t\to\infty} \sum_{n=2}^{\infty} b_ne^{-r_nt} = 0,
$$
which implies $b_1=0$.  Continuing the argument, we can reach 
$b_n=0$ for all $n\in \mathds{N}$. \hfill$\blacksquare$

\bigskip
\noindent
\textbf{Proof of Theorem~\ref{th.uniq2}:} By \eqref{eq.uon}, we see that 
$$
\sum_{n=1}^{\infty} \frac{4(1-a)^2\vert U_n^0(a)\vert^2}
{j_n^4 J_0'(j_n)^2}
= \sum_{n=1}^{\infty} \vert u^0_n\vert^2 < \infty,
$$
that is,

$$
\sum_{n=1}^{\infty} \left\vert \frac{U_n^0(a)}
{j_n J_0'(j_n)}\right\vert^2\frac{1}{n^2}
\le C\sum_{n=1}^{\infty} \left\vert \frac{U_n^0(a)}
{j_n J_0'(j_n)}\right\vert^2\frac{1}{j_n^2} < \infty.
$$

\noindent
Therefore, the assumption of Lemma~\ref{lemma5.1} with $\gamma=1$ is satisfied.
Hence, the time-analyticity yields 

$$
\sum_{n=1}^{\infty} \frac{U_n^0(a_1)}{j_n J_0'(j_n)} e^{-\lambda_nt}
- 
\sum_{n=1}^{\infty} \frac{\widetilde{U}_n^0(a_2)}{j_n J_0'(j_n)} e^{-\mu_nt}
= 0, \quad t>t_1.
$$
First assume that $\{ \lambda_n\}_{n\in  \mathds{N}} \cap \{ \mu_n\}_{n\in  \mathds{N}}
= \emptyset$, where $\lambda_n$ and $\mu_n$ are given by  \eqref{lambdamu}. Then, we can renumber $\lambda_n, \mu_n$ as 
$0 < \theta_1 < \theta_2 < \cdots$, and 
$$
\sum_{m=1}^{\infty} V_me^{-\theta_mt} = 0, \quad t>t_1,
$$
where $V_m$ is given by 
$\frac{U_n^0(a_1)}{j_n J_0'(j_n)}$ or $\frac{\widetilde{U}_n^0(a_2)}{j_n J_0'(j_n)}$.

Applying Lemma~\ref{lemma5.1}, we see $V_m=0$ for all $m\in  \mathds{N}$. 
This yields $u^{a_1} \equiv 0$ and $u^{a_2} \equiv 0$ in $(0,1) \times 
(0,\infty)$.  Since we assume that $u^{a_1} \not\equiv 0$ or 
$u^{a_2} \not\equiv 0$ in $(0,1) \times (0,\infty)$, this is impossible.
Thus, we verified that $\{ \lambda_n\}_{n\in  \mathds{N}} \cap \{ \mu_n\}
_{n\in  \mathds{N}} \ne \emptyset$, that is, we can find $m_1, m_2 \in \mathds{N}$ such that 
$\lambda_{m_1} = \mu_{m_2}$, and we can complete the proof of Theorem~\ref{th.uniq2}.

\hfill$\blacksquare$

\section{Numerical results}\label{sec:numerics}

In this section we will show some numerical results related to the identification of a degeneracy point $a\in(0,1)$ in \eqref{eq.cp}. 

More precisely, given $T>0$ and $\beta(t)$, we will present some numerical tests for different initial data $w_0$, so as we can find $a\in (0,1)$ such that the solution to \eqref{eq.cp} for some $t_0\in (0,T)$ satisfies  
\begin{equation}\label{obs1xo}
w_x(1,t_0) = \beta(t_0).
\end{equation}

In order to reconstruct $a$, we will reformulate the Inverse Problem~\ref{IP_refor} as an optimization problem. With fixed small $\delta>0$, let us consider the admissible set 

\begin{equation}\label{uad}
\mathcal{U}^{a}_{ad} =\{ a: a\in (\delta, 1-\delta) \} 
\end{equation}
and a functional $J : a\in \mathcal{U}^{a}_{ad} \mapsto \mathds{R}$ given by

\begin{equation}\label{J1x0}
J(a) = \dfrac{1}{2}\displaystyle |\beta(t_0) - w^a_x(1,t_0)|^2
\end{equation}
for some $t_0\in (0,T)$.
\noindent
The related optimization problem is the following:

\begin{equation}\label{optpb1x0}
\left\{\begin{array}{l}
\text{Minimize  $J(a) $ }   \\[1mm]
\text{where $a\in \mathcal{U}^{a}_{ad} $ and $w^a$ satisfies \eqref{eq.cp}. }
\end{array}\right.
\end{equation}

\medskip
Let us present several numerical tests that illustrate the theoretical results from the previous sections.  
The \texttt{fmincon} function from MATLAB Optimization ToolBox (the gradient method) will be used in order to solve the constrained optimization problem~\eqref{optpb1x0}. 

\begin{figure}[h!]
\begin{minipage}[t]{0.49\linewidth}
\noindent
\includegraphics[width=0.98\linewidth]{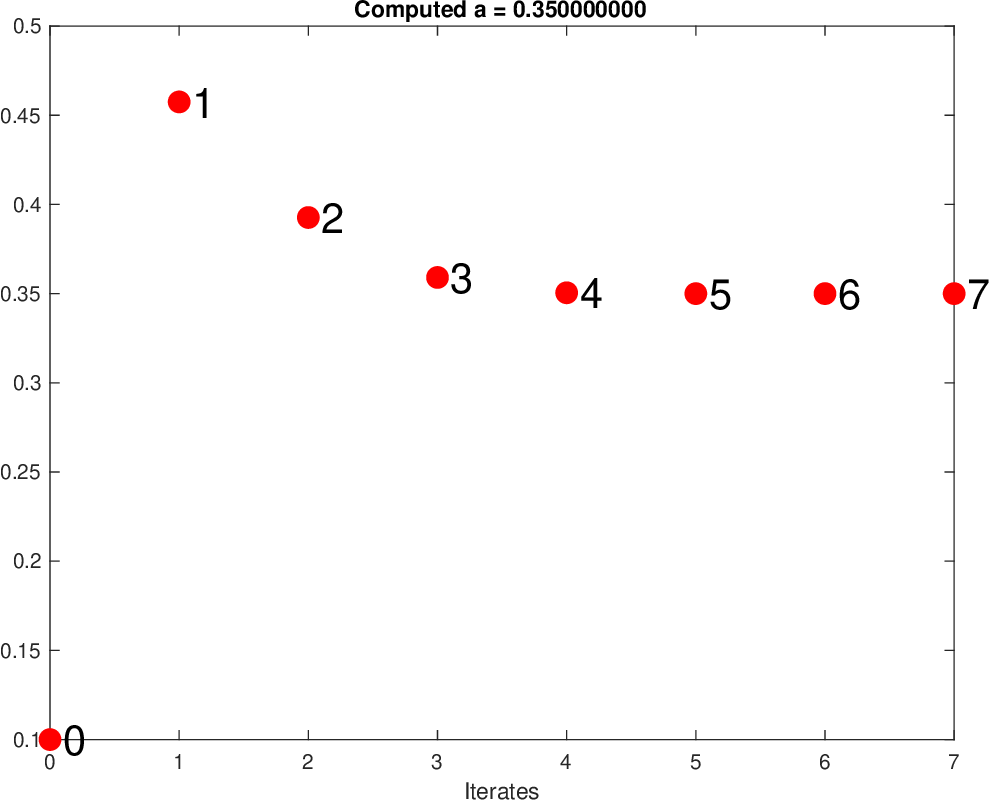}
\caption{Test~\ref{test1}, $w_0=1$. Iterations in the computation of $a$ by \texttt{trust-region-reflective} algorithm, $aini=0.1$. }
\label{a_iter_u01}
\end{minipage}
\hfill
\begin{minipage}[t]{0.49\linewidth}
\noindent
\includegraphics[width=0.98\linewidth]{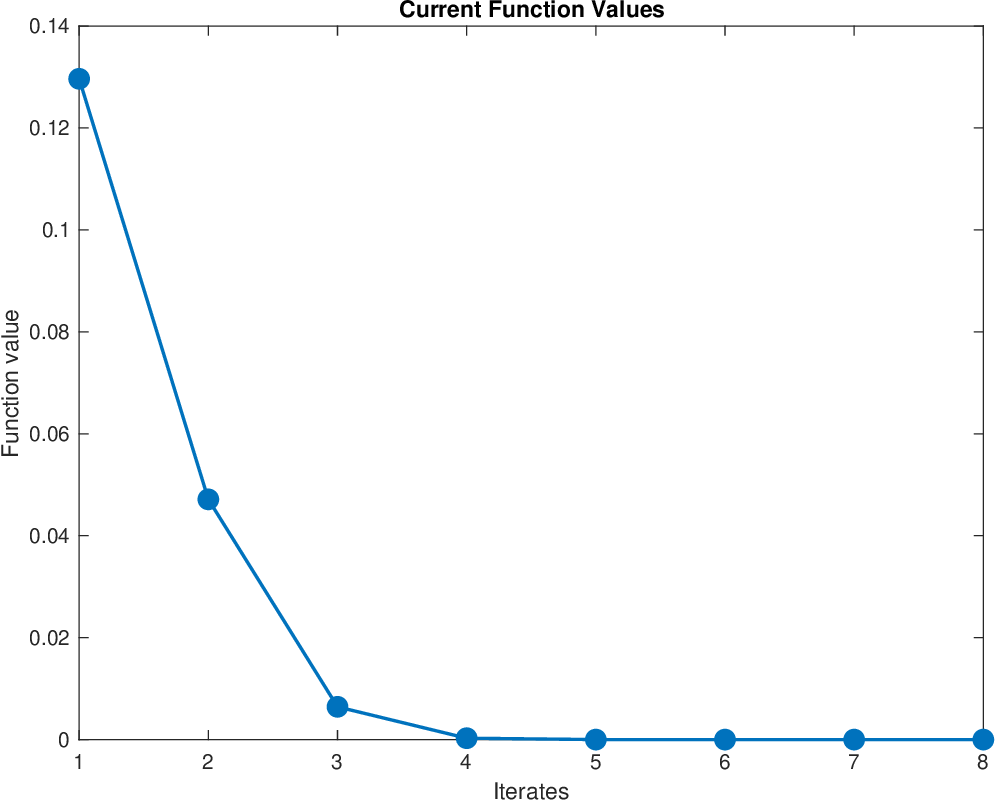}
\caption{Test~\ref{test1}, $w_0=1$. Evolution of the cost in \texttt{trust-region-reflective} algorithm, $aini=0.1$.}
\label{fvalues_a_u01}
\end{minipage}
\end{figure}

\begin{test} \label{test1}
The goal is to reconstruct the  degeneracy point $a$ for the initial data for the initial data $w_0=1$.
We will take $aini=0.1$ and also $aini=0.6$ as initial guesses for the recovering the desired value of $a_d=0.35$ by the minimization algorithm and $T=10$.  We will take in \eqref{J1x0} $t_0 = 0.05$. Notice that in this case, due to Proposition~\ref{prop.case1}, the time for the observation $t$ in \eqref{J1x0} can be arbitrarily small. 

The  numerical results can be seen in Figures~\ref{a_iter_u01}, \ref{fvalues_a_u01}. The round points correspond to iterations during the optimization algorithm.  With solid line, we have represented the evolution of the cost.  
\end{test}
\medskip

\begin{table}[h!]
\centering
\renewcommand{\arraystretch}{1.2}
\caption{Some results for different initial guess and algorithms, Tests~\ref{test1}.}
\medskip
{\begin{tabular}{cccccc } \hline
Initial guess    & Computed $a$ & Algorithm & Iterates &  Cost 
\\
\hline
0.1                 & 0.3503244447468097 & interior-point                     & 6 &  1.e-7    \\
\hline
0.1                &  0.3499999999999553       &trust-region-reflective   &  7   &  1.e-27     \\
\hline
0.1                &  0.3500000021587790       & active-set                    & 8    &  1.e-17    \\
\hline\hline
0.6                &   0.3500000000000531      & interior-point               & 15  &  1.e-30    \\
\hline
0.6                &  0.3500000000001043       & trust-region-reflective   & 7  &  1.e-26    \\
\hline
0.6               & 0.3499999911140417         & active-set                       & 7  &  1.e-16  \\
\hline
\end{tabular}}
\label{Table1}
\end{table}


\bigskip

   In Table~\ref{Table2} we can see the evolution of the cost when we introduce random noises in the target. These results correspond to the \texttt{trust-region-reflective} algorithm.

\medskip
\begin{table}[h!]
\centering
\renewcommand{\arraystretch}{1.2}
\caption{Evolution of the cost with random noises in the target, Test~\ref{test1} with  $aini=0.1$.}
\medskip
{\begin{tabular}{cccc} \hline
random noise  & Cost in \texttt{trust-region-reflective}  & Iterates &  Computed $a$ \\
\hline
1\%         &  1.e-24          &    7      & 0.339169874402234 \\ \hline
0.1\%      &  1.e-25           &   7      & 0.3488748208839999\\ \hline
0.01\%      &  1.e-27          &   7     & 0.3500145020619460   \\\hline
0.001\%    &  1.e-27         &   7       & 0.3500141069331907  \\ \hline
0\%           &   1.e-27        &    7      & 0.3499999999999553   \\  \hline
\end{tabular}}
\label{Table2}
\end{table}

\begin{figure}[h!]
\begin{minipage}[t]{0.49\linewidth}
\noindent
\includegraphics[width=0.98\linewidth]{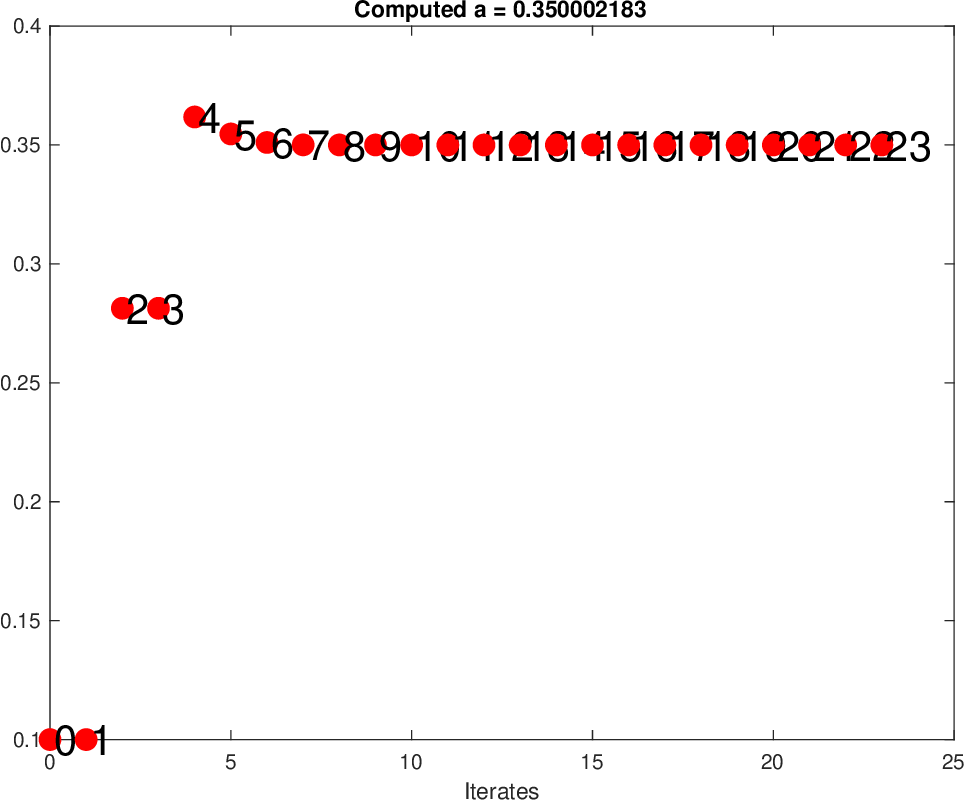}
\caption{Test~\ref{test2}, $w_0(x)=1-x$. Iterations in the computation of $a$ in \texttt{trust-region-reflective} algorithm. }
\label{a_iter_u01mx}
\end{minipage}
\hfill
\begin{minipage}[t]{0.49\linewidth}
\noindent
\includegraphics[width=0.98\linewidth]{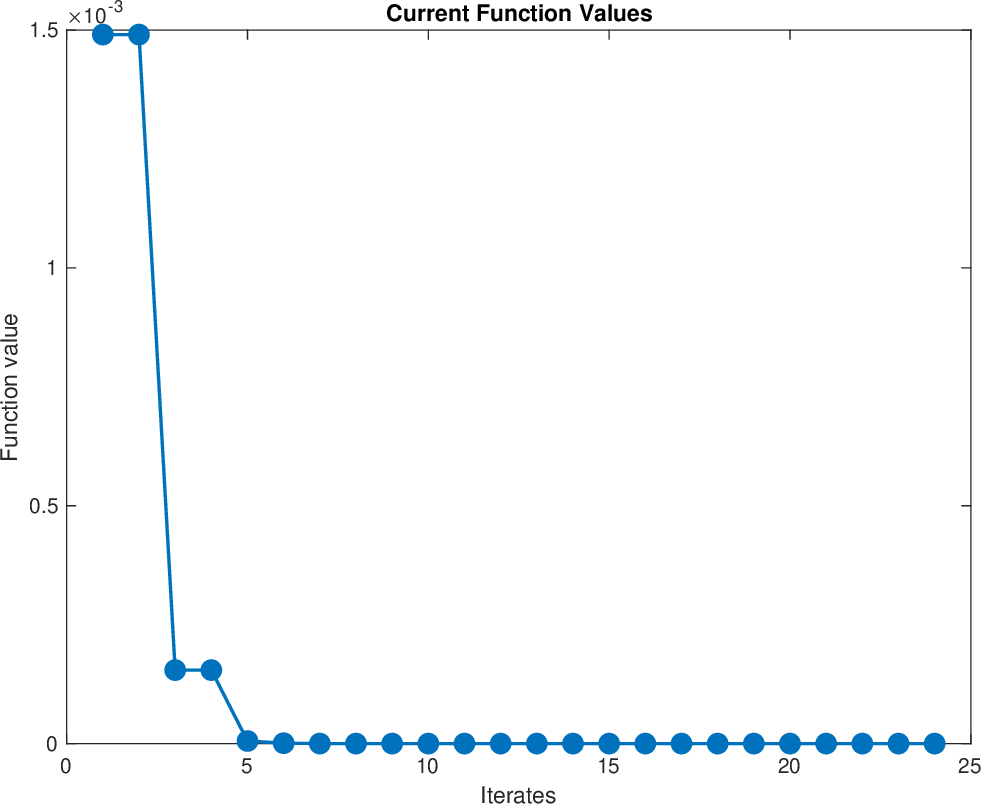}
\caption{Test~\ref{test2}, $w_0(x)=1-x$. Evolution of the cost in \texttt{trust-region-reflective} algorithm.}
\label{fvalues_a_u01mx}
\end{minipage}
\end{figure}

\bigskip
\begin{test} \label{test2}
In this test we reconstruct the degeneracy point $a$ for the initial data $w_0(x)=1-x$ from the second theoretical example. We will take $aini=0.1$ as an initial guesses for the recovering the desired value of $a_d=0.35$ by the minimization algorithm.  Again, we have taken in \eqref{J1x0} $t_0 = 0.05$.

The  numerical results obtained by application of the \texttt{trust-region-reflective} algorithm can be seen in Figures~\ref{a_iter_u01mx}, \ref{fvalues_a_u01mx} and Table~\ref{Table3}. Again, the round points correspond to iterations of the optimization algorithm and the solid line the evolution of the cost.  
\end{test}
\medskip

\medskip
\begin{table}[h!]
\centering
\renewcommand{\arraystretch}{1.2}
\caption{Evolution of the cost with random noises in the target, Test~\ref{test2}.}
\medskip
{\begin{tabular}{cccc} \hline
\% random noise               & Cost  & Iterates &  Computed $a$ \\
\hline
1\%         &  1.e-32            &    10      & 0.3360403710906295 \\ \hline
0.1\%      &  1.e-28            &   11      & 0.3510948436764733\\ \hline
0.01\%      &  1.e-27         &   10        & 0.3498074191585610   \\\hline
0.001\%    &  1.e-15        &    27       & 0.3360403710906295  \\ \hline
0\%           &   1.e-14       &    23        & 0.3500021829579620   \\  \hline
\end{tabular}}
\label{Table3}
\end{table}

\begin{figure}[h!]
\centering
\includegraphics[width=0.58\linewidth]{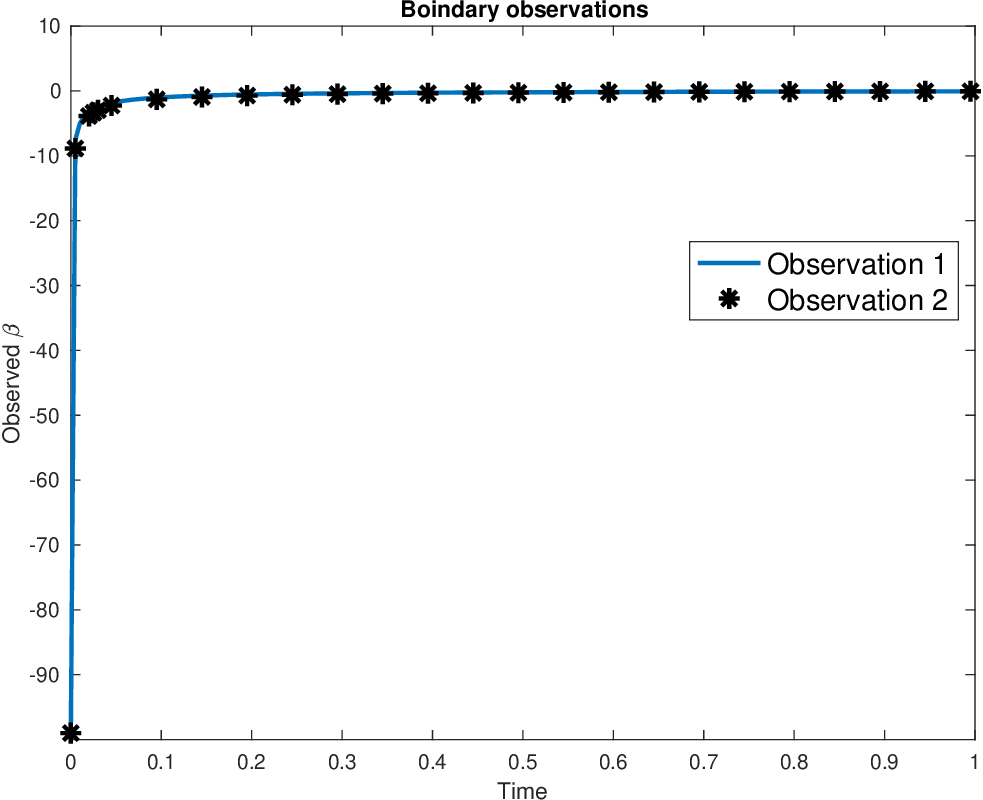}
\caption{$w_0(x)=x$, observations of the values $w^a_x(1,t_0)$, Test~\ref{test3}.}
\label{fig.BoundaryObs}
\end{figure}

\begin{test}\label{test3}

Let us consider $w_0(x) = x$. In this cases, in general we cannot ensure the uniqueness of the inverse problem (see Figure~\ref{fig.normderux1}).   

In order to illustrate this issue, we will consider two cases: 

\begin{description}
\item[Case 1:] Starting the optimization algorithm from $aini = 0.1$ with the goal to recover the desired value $a_d = 0.163$, after the application of the \texttt{trust-region-reflective} algorithm, we obtain the computed value $a = 0.1629999999999927$ with the cost 1.e-27. 
\item[Case 2:] Starting the optimization algorithm from $aini = 0.4$ with the goal to recover the desired value $a_d = 0.379$, after the application of the \texttt{trust-region-reflective} algorithm, we obtain the computed value $a = 0.3789957650933454$ with the cost 1.e-11. 
\end{description} 

However, Figure~\ref{fig.BoundaryObs}  indicates that different values of $a$ produce the same boundary observations $\beta(t_0) = w^a_x(1,t_0)$. 
\end{test}
\medskip

\bigskip

\appendix 
\section{Proof of Lemma~\ref{lemma.Bessel}}\label{sec.appendix}

For completeness, we present  the proof of Lemma~\ref{lemma.Bessel}. 

\medskip
\noindent
\textit{a)} For all $n = 1, 2\dots $, we can write

\begin{equation}\label{eq.Besselder}
\begin{split}
\dfrac{d}{dz}\big (z^n J_n(z)\big) & = \displaystyle\sum_{k=0}^\infty \dfrac{(-1)^k}{k! (n+k)!} \frac{2n+2k}{2^{n+2k}}z^{2n+2k-1}
\\[3mm] 
 & = 
 \displaystyle\sum_{k=0}^\infty z^n \dfrac{(-1)^k}{k! (n+k-1)!} \Big(\frac{z}{2}\Big)^{n+2k-1} : = z^n J_{n-1}(z).
 \end{split}
\end{equation}
Consequently, for all $x>0$ we obtain
\begin{equation}\label{eq.Besselint}
\displaystyle  \int_0^x s J_0(s) \, ds = s J_1(s) \Big|_{0}^x =x J_1(x).
\end{equation}
%

\medskip
\noindent
\textit{b)} Similarly, we obtain 

\begin{equation}\label{eq.Besselder2}
\dfrac{d}{dz}\big (z^{-n} J_n(z)\big)  = -z^{-n} J_{n+1}(z).
\end{equation}

\noindent 
Therefore, we have
\[
n z^{n-1} J_n(z) + z^n J_n'(z) = z^n J_{n-1}(z). 
\]
Thus, we have 
\[
\begin{cases}
J_n'(z) +\dfrac{n}{z} J_n(z) = J_{n-1}(z),
\\[4mm] 
J_n'(z) -\dfrac{n}{z} J_n(z) = J_{n+1}(z),
\end{cases}
\]
which in turn implies 
\[
\begin{cases}
J_{n-1}(z) + J_{n+1}(z) =\dfrac{2n}{z} J_n(z),
\\[3mm] 
J_{n-1}(z) - J_{n+1}(z) =2J'_n(z).
\end{cases}
\]

\medskip
\noindent 
\textit{c)} From \eqref{eq.besselJ0}, we have 

\[
J'_0(z) = \displaystyle\sum_{k=1}^\infty \dfrac{(-1)^k}{k! (k-1)!} \dfrac{2 z^{2k-1}}{2^{2k}} 
=
- \displaystyle\sum_{k=0}^\infty \dfrac{(-1)^k}{k! (k+1)!} \Big(\dfrac{z}{2}\Big)^{1+2k}  := -J_1(z). 
\]

\medskip
\noindent 
\textit{d)} Taking into account, the property \textit{c)}, we have 
\[
\begin{split}
\displaystyle\int_0^{j_n} s J_0^2(s) \, ds & = \Big[\dfrac{s^2}{2}J_0(s)^2  \Big]_{0}^{j_n}
- \displaystyle\int_0^{j_n} s^2 J_0(s) J_0'(s)\, ds 
\\[3mm]
 & = 
 \displaystyle\int_0^{j_n} (sJ_1(s))' s J_1(s)(s) \, ds 
 = \Big[\dfrac{s^2}{2}J_1(s)^ 2 \Big]_{0}^{j_n}.
\end{split}
\]

Therefore, we see \textit{d)}. 

\medskip
\noindent 
\textit{e)} Using $\displaystyle\int_0^{x} s J_0(s) \, ds =\displaystyle\int_0^{x} \dfrac{d}{ds}(s J_1(s)) \, ds$ we have  \textit{e)}.

\bigskip
\noindent 
\textit{f)}  The properties of zeros of the function $J_\nu$ can be found in~\cite{Lebedev} and the bounds for zeros of $J_0$ are given in~\cite{LorchMuldoon}, see also~\cite{CMV2}.

\bigskip
\noindent 
\textit{g)} It is know (see~\cite{Lebedev}) that 
$$
J_0(z) \sim \sqrt{\frac{2}{\pi z}} \cos \left( z - \frac{\pi}{4}\right)
$$
as $\vert z\vert \to \infty$.  Hence,
$$
\vert J_0(s)\vert \le C_1\left( \frac{1}{\sqrt{s}} + 1\right), \quad 
s > 0.
$$
Therefore Lemma~\ref{lemma.Bessel} \textit{d)} and \textit{f)} implies 
$$
\frac{j_n^2}{2}J_0'(j_n)^2 = \int^{j_n}_0 sJ_0^2(s) ds
\le \int^{j_n}_0 sC_2\left( 1 + \frac{1}{s}\right) ds 
\le C_3(j_n^2+j_n),
$$
where $C_i$ are positive constants. 

Hence, $J_0'(j_n)^2 \le 2C_3\left( 1 + \frac{1}{j_n}\right)$ and 
we see the property \textit{g)}.

\hfill $\blacksquare$

\section*{Acknowledgments}
The first author was supported in part by the National Group for Mathematical Analysis, Probability and Applications (GNAMPA) of the Italian Istituto Nazionale di Alta Matematica ``Francesco Severi'' and by the Excellence Department Project awarded to the Department of Mathematics, University of Rome Tor Vergata, CUP E83C23000330006. The second author was partially supported by PROYECTO PID2020-114976GB-I00, PLAN ESTATAL 2021 PROYECTO PY20\_01125, PAIDI 2021. The third author is supported by Grant-in-Aid for Scientific Research (A) 20H00117 and Grant-in-Aid for Challenging Research (Pioneering) 21K18142 from Japan Society for the Promotion of Science.

\bigskip

\end{document}